\documentclass[11pt]{article}
\usepackage[T2A]{fontenc}
\usepackage[cp1251]{inputenc}
\usepackage{amsmath}
\usepackage{amssymb}
\usepackage{amsthm}
\usepackage{arydshln}
\def\a{\mathbf{a}}
\def\hh{\mathbf{h}}
\def\pp{\mathbf{p}}

\def\m{\mathbf{m}}
\def\f{\mathbf{f}}

\def\b{\mathbf{b}}
\def\N{\mathbb{N}}
\def\Z{\mathbb{Z}}
\def\R{\mathbb{R}}
\def\C{\mathbb{C}}

\def\etta{\boldsymbol{\eta}}
\def\zetta{\boldsymbol{\zeta}}
\def\llambda{\boldsymbol{\lambda}}
\def\ggamma{\boldsymbol{\gamma}}
\def\aalpha{\boldsymbol{\alpha}}
\def\bbeta{\boldsymbol{\beta}}
\def\llambda{\boldsymbol{\lambda}}

\def\ddelta{\boldsymbol{\delta}}
\def\rrho{\boldsymbol{\rho}}
\newtheorem{Lemma}{\hspace*{\parindent}Lemma}
\newtheorem{Theorem}{\hspace*{\parindent}Theorem}

\title{Beyond the beta integral method: transformation formulas for hypergeometric functions via Meijer's $G$ function}
\author{D.B.\:Karp$^{\rm a,b}$\footnote{Corresponding author. E-mail: D. Karp -- \emph{dimkrp@gmail.com},
E.\:Prilepkina --  \emph{pril-elena@yandex.ru}}~~and
E.G.\:Prilepkina$^{\rm b,c}$
\\[10pt]
\\
\small{\textit{$\phantom{1}^a$Holon Institute of Technology, Holon, Israel}}
\\
\small{\textit{$\phantom{1}^b$Far Eastern Federal University, Vladivostok, Russia}}
\\
\small{\textit{$\phantom{1}^c$Institute of Applied Mathematics,
FEBRAS,  Vladivostok,   Russia}}}
\date{}

\begin{document}

\maketitle

\begin{center}
\parbox{12cm}{
\small\textbf{Abstract.}
The beta integral method proved itself as a simple 	nonetheless powerful method of generating hypergeometric identities at a fixed argument.  
In this paper we propose a generalization by substituting the beta density with a particular type of Meijer's $G$ function.  
By application of our method to known transformation formulas we derive about forty hypergeometric identities, 
majority of which are believed to be new. We further apply some of these transformations to obtain some new summation formulas. }
\end{center}

\bigskip


\bigskip


\bigskip
\section{Introduction}

Summation and transformation formulas for hypergeometric functions at a fixed argument are important in combinatorics  \cite{AndrStan,Chu2002,Haglund}, analysis \cite{ChoChungYun,Koepf2007,Skwarczynski}, physics \cite{Minton,RaoBook,RDN}, computer science \cite{GreeneKnuth} and many other fields \cite{AAR,Seaborn}.
As most summation formulas are particular or limiting cases of some transformation formulas, the latter turn out to be of a higher significance.  The main developments until the end of 1930ies were summarized by W.N.\:Bailey in the fundamental monograph  \cite{Bailey}.  His student, Lucy Joan Slater attributes to L.J.\:Rogers the statement that after Bailey's work ''nothing remains to be done in the field of hypergeometric series'' \cite[p.40]{Slater}.  In his work, Bailey gave a number of methods for deriving and proving such transformation formulas, including series rearrangements, contour integrals, equating coefficients in an identity involving free argument, the ''Bailey method'' \cite[Lemma~3.4.2]{AAR} and the Bailey chains, \cite[Chapter~12]{AAR}. Later on, an important extension to this toolbox was provided by the algorithms of symbolic computation \cite{Koepf,Koepf2007},\cite[section~3.11]{AAR}, techniques based on Lagrange inversion theorem \cite{GesselStant} and  Abel's lemma \cite{ChuWang2009} and various other methods \cite{Chen,Choi}.  Another simple, but frequently very effective method for obtaining transformation formulas at a fixed argument from an identity involving a free argument consists in integrating such identity with respect to the beta density.  It pops up in the literature on various occasions but was fully automated and systematically applied by Krattenthaler and Rao in \cite{KrRao2003} and was given the name ''the beta integral method'' by these authors \cite[Chapter~8]{RaoBook}.  The main idea of this work is to generalize this method by substituting the beta density with a density expressed in terms of Meijer's $G$ function, of which the beta density is a particular case.  Unlike the beta integral method, however, this approach does not automatically lead to a transformation formula. The reason behind this phenomenon is that for the beta integral method to work one only needs ${}_2F_{1}(1)$ to be summable in terms of gamma functions which is always the case by the celebrated Gauss formula. In contrast, for an application of $G$ function integral method proposed here one needs a summation formula for the generalized hypergeometric function ${}_{p+1}F_{p}(1)$ with $p\ge2$, which imposes severe parameter restrictions.  These restrictions in many cases contradict the parameter structure dictated by the $G$ function integral method. There is a number of cases, however, when these two requirements are compatible and we are led to transformation and summation formulas for the generalized hypergeometric functions evaluated at a specific value of the argument (typically at unity).

It is convenient to introduce an extended definition of the  hypergeometric series by
\begin{equation}\label{eq:pFqdefined}
F\left.\!\left(\!\begin{array}{l}\a \\\b \end{array}{\vdots}\,  P\, \right\vert x\right)
=\sum\limits_{n=0}^{\infty}\frac{(a_1)_n(a_2)_n\cdots(a_{p})_n}{(b_1)_n(b_2)_n\cdots(b_q)_nn!}P(n)x^n,
\end{equation}
where $\a=\{a_1,\ldots,a_p\}$, $\b=\{b_1,\ldots,b_q\}$ are complex parameter vectors such that $b_j$ never equals a non-positive integer and $P(n)$ could be any function of $n$, but in this paper it will always be a polynomial of a fixed degree $m$. The expression $(a)_n=\Gamma(a+n)/\Gamma(a)$ is  the standard Pochhammer's symbol (or rising factorial).   In this case
it is straightforward to check that
$$
P(n)=P(0)\frac{(1-\llambda)_{n}}{(-\llambda)_{n}},
$$
where $\llambda=(\lambda_1,\ldots,\lambda_m)$ is the vector of zeros of the polynomial $P$ (repeated if necessary according to the multiplicity) and the shorthand notation for the product $(-\llambda)_{n}=(-\lambda_1)_n(-\lambda_2)_n\cdots(-\lambda_m)_n$ is used here and henceforth. Hence, \eqref{eq:pFqdefined} can be rewritten as 
\begin{equation}\label{eq:FP-roots}
F\left.\!\left(\!\begin{array}{l}\a \\\b \end{array}{\vdots}\, P\, \right\vert x\right)
=P(0){}_{p+m}F_{q+m}\left.\!\left(\!\begin{array}{l}\a, 1-\llambda \\\b, -\llambda \end{array}\right\vert x\right)
\end{equation}
- a generalized hypergeometric function with $m$ unit shifts in the parameters.  We will use both ways of writing $F$ interchangeably.
This extended definition has been recently employed by Maier \cite{Maier2019} and is equivalent to the concept of ''hypergeometrization'' introduced  a bit earlier by Blaschke \cite{Blaschke}.  We also found it convenient to omit the indices of the hypergeometric functions, as the dimensions of the parameter vectors are usually clear from the context. However, we will use the traditional notation ${}_pF_{q}$ when dealing with specific numerical values of $p$ and $q$ to make the formulas more accessible to a reader not interested in further details.
To avoid poles in the denominators we will always assume that $b_1,\ldots,b_q$ do not take  non-positive integer values.
Finally, omitted argument of the generalized hypergeometric function signifies the unit argument throughout the paper.

The paper is organized as follows.  In the next section we describe the general framework and present a list of summation and transformation formulas that will serve as a raw material for our machinery.  In section~3 we present the transformation formulas obtained by the $G$ function integral method applied to the identities presented in subsection~2.3. We group the formulas in accordance with the values of the parameters $u$ and $v$ in \eqref{eq:generaltrans}.
We included both the formulas we could not locate in the literature and a few well-known transformations to illustrate the power of the method. We further added a reference each time we were aware of it.   It is typically rather difficult to claim that a hypergeometric transformation is new, as the literature is vast and there could always be a hidden trick how a ''new'' transformation can be derived  from a known one.  Hence, we simply present all formulas that we reckoned as interesting with the hope that some one them are indeed new.

\section{$G$ function integral method: preparation}

\subsection{General description of the method}
We will use the standard symbols $\N$, $\Z$, $\R$ and $\C$ to denote the sets of natural, integer, real, and complex numbers, respectively. Similarly to the beta integral method, we will start with a transformation formula of the form
\begin{equation}\label{eq:generaltrans}
F\left.\!\left(\begin{matrix}\aalpha\\\bbeta\end{matrix}\:\right\vert
Mx^{w}\right)
=(1-x)^{\lambda}F\left.\!\left(\begin{matrix}\ddelta\\\ggamma\end{matrix}\:\right\vert
\frac{Dx^{u}}{(1-x)^v}\right)
\end{equation}
valid for $0<x<1$.  Here $\ddelta$, $\ggamma$ and $\lambda$ are functions of $\aalpha$, $\bbeta$; $w,u\in\N$, $v\in\Z$, $M,D$ are constants.  A list of transformations of the form \eqref{eq:generaltrans} will be given in the following subsection.  Here we just note that the cases
$$
(w,u,v)\in\{(1,1,0), (1,1,1)\}
$$
correspond to the Euler-Pfaff linear-fractional transformations and their generalizations to hypergeometric functions with integral parameter differences \cite{KPChapter2019,MP2013}.  The quadratic transformations include the cases
$$
(w,u,v)\in\{ (1,1,2), (2,1,1),  (2,1,2), (1,2,2), (1,1,-1), (1,2,1)\}.
$$
Some cubic transformations  \cite{Askey1994} also have the form \eqref{eq:generaltrans}. These cases have been explored by us in \cite{CKP-Lob}.

The beta integral method consists in multiplication of the transformation formula \eqref{eq:generaltrans} by the beta density $x^{d-1}(1-x)^{e-1}$ and term-wise integration from $0$ to $1$. In this work we substitute the beta density by the Meijer-N{\o}rlund function $G^{p,0}_{p,p}$ of which it is a particular $p=1$ case. This function is defined by the Mellin-Barnes integral of the form
\begin{equation}\label{eq:G-defined}
G^{p,0}_{p,p}\!\left(\!z~\vline\begin{array}{l}\b\\\a\end{array}\!\!\right)\!\!:=
\\
\frac{1}{2\pi{i}}
\int\limits_{\mathcal{L}}\!\!\frac{\Gamma(\a\!+\!s)}{\Gamma(\b+\!s)}z^{-s}ds.
\end{equation}
The shorthand notation $\Gamma(\a\!+\!s)=\prod_{j=1}^{p}\Gamma(a_{j}+s)$ is used here and henceforth. Details regarding the choice of the contour ${\mathcal{L}}$ can be found in many standard reference books \cite[section~5.2]{LukeBook}, \cite[16.17]{NIST}, \cite[8.2]{PBM3} and our papers \cite{KLJAT2017,KPSIGMA}, which also contain a list of properties of $G^{p,0}_{p,p}$.  In particular,  to perform the term-wise integration we will need the integral evaluation \cite[p.50]{KLJAT2017}
\begin{equation}\label{eq:Gpowerint}
\frac{\Gamma(\b)}{\Gamma(\a)}\int_0^1x^{\nu}(1-x)^{\mu}G^{p,0}_{p,p}\!\left(\!x\left|\!\begin{array}{l}\b-1\\\a-1\end{array}\right.\!\!\!\right)dx
=\frac{(\a)_{\nu}}{(\b)_{\nu}}{}_{p+1}F_{p}\left.\!\!\left(\!\begin{matrix}-\mu,\a+\nu\\\b+\nu\end{matrix}\right.\right),
\end{equation}
where for any $\nu$ the Pochhammer's symbol is given by $(a)_{\nu}=\Gamma(a+\nu)/\Gamma(a)$ and $(\a)_{\nu}$ is the shorthand notation for the product $\prod_{j=1}^{p}(a_j)_{\nu}$. The above formula is true if $\Re(\a+\nu)>0$ and $\Re(s(\a,\b)+\mu)>0$ (understood element-wise), where $s(\a,\b)$ here and below signifies \emph{the parametric excess}
\begin{equation}\label{eq:parametricexcess}
s(\a,\b)=\sum_{j=1}^{p}(b_{j}-a_{j}).
\end{equation}

An application of these ideas leads to the following ''master lemma''.
\begin{Lemma}\label{lm:master1}
Assume that \eqref{eq:generaltrans} holds for  $x\in(0,1)$. Suppose further that $\ddelta$ or $\a$ contain a negative integer or $v=0$, $D=1$, and
\begin{equation}\label{eq:v0condition}
\Re(\a)>0~\&~\Re(s(\a,\b)+\lambda)>0~\&~\Re(s(\a,\b)+s(\ggamma,\ddelta)+\lambda)>0.
\end{equation}
Then
\begin{equation}\label{eq:master1}
F\left.\!\!\left(\!\begin{matrix}\aalpha,\Delta(\a,w)\\\bbeta,\Delta(\b,w)\end{matrix}\right\vert M\right)
=\sum\limits_{k=0}^{\infty}\frac{(\ddelta)_k(\a)_{uk}D^k}{(\ggamma)_k(\b)_{uk}k!}
F\left.\!\!\left(\!\begin{matrix}-\lambda+vk,\a+uk\\\b+uk\end{matrix}\right.\right),
\end{equation}
where $\Delta(a,w)=(a/w,a/w+1/w,\ldots,a/w+(w-1)/w)$.
\end{Lemma}

\textbf{Remark.} Only the last restriction in \eqref{eq:v0condition} is required for convergence of the series on the right hand side of \eqref{eq:master1}, while the first two restrictions are needed for the derivation only and can  typically be removed by analytic continuation once a transformation formula has been obtained.

\begin{proof} For the proof multiply \eqref{eq:generaltrans} by
$$
\frac{\Gamma(\b)}{\Gamma(\a)}G^{p,0}_{p,p}\!\left(\!x\left|\!\begin{array}{l}\b-1\\\a-1\end{array}\right.\!\!\!\right)
$$
and integrate both sides term-wise from $0$ to $1$ using \eqref{eq:Gpowerint} with $\mu=0$ on the left hand side and $\mu=-\lambda-vk$ on the right hand side.
Further, apply
\begin{equation}\label{eq:multiple-poch}
(a)_{lk}=l^{lk}(a/l)_{k}((a+1)/l)_{k}\cdots((a+l-1)/l)_{k}=l^{lk}\Delta(a,l)_{k},
\end{equation}
valid for each $l\in\N$, on the RHS. If $\ddelta$ or $\a$ contain a negative integer, the summation terminates and term-wise integration is permitted.

Otherwise, if $v=0$, $D=1$ define
$$
F(k):=\frac{\Gamma(\a+uk)}{\Gamma(\b+uk)}{}_{p+1}F_{p}\left.\!\!\left(\!\begin{matrix}-\lambda,\a+uk\\\b+uk\end{matrix}\right.\right).
$$
By the change of variable $x=e^{-t}$ in \eqref{eq:Gpowerint} we obtain
$$
F(k)=\int\limits_{0}^{\infty}e^{-tuk}G^{p,0}_{p,p}\!\left(\!e^{-t}~\vline\begin{array}{l}\b-1\\\a-1\end{array}\!\!\right)
(1-e^{-t})^{\lambda}e^{-t}dt=\int\limits_{0}^{\infty}e^{-kf(t)}g(t)dt,
$$
where
$$
f(t)=tu,~~~g(t)=(1-e^{-t})^{\lambda}e^{-t}G^{p,0}_{p,p}\!\left(\!e^{-t}~\vline\begin{array}{l}\b-1\\\a-1\end{array}\!\!\right)
=t^{\lambda+s(\a,\b)-1}\hat{g}(t),
$$
and the function $\hat{g}(t)$ is analytic near $t=0$ with $\hat{g}(0)\ne0$ according to \cite[(11)]{KLJAT2017}. First two conditions in \eqref{eq:v0condition} make sure that the integral converges. The function $f(t)$ has the minimum at $t=0$ with $f(0)=0$ and $f'(0)\ne0$. An application of Watson's lemma \cite[Theorem~C.3.1]{AAR} then yields:
$$
F(k)\sim C_1k^{-s(\a,\b)-\lambda}~\text{as}~k\to\infty.
$$
In view of
$$
\frac{(\ddelta)_k}{(\ggamma)_kk!} \sim C_2 k^{-s(\ggamma,\ddelta)-1},
$$
the series on the right hand  side of \eqref{eq:master1} absolutely converges if the third condition in  \eqref{eq:v0condition} is satisfied, so that term-wise integration is justified. \end{proof}

\subsection{Summation formulas}

In this section we will list the cases when the hypergeometric function of the right hand side of \eqref{eq:master1} is summable in terms of gamma functions.  These cases hinge on the classical summation theorems, their extensions, and the following lemma for hypergeometric functions with integral parameter differences (IPD type).  A related formula can be found in our paper \cite[Theorem~3.2]{KPITSF2018}. Both in this  lemma and in the sequel we will use the notation $\pp=(p_1,\ldots,p_l)\in\N^l$,  $p=p_1+\cdots+p_{l}$ and $\hh\in\C^{l}$.  Let us emphasize that all formulas presented in this section are essentially known results rewritten in the form convenient for further application in Section~3 which is devoted to new results. 
\begin{Lemma}\label{lm:IPDsummation}
Suppose $l\in\N$,  $u$, $v$ are integers.
Then for $k\in\N$ such that $\Re(e+\lambda-d-p-vk)>0$ or if hypergeometric function $F$ terminates, we have
\begin{equation}\label{eq:IPDsummatiom}
F\!\left(\begin{matrix}-\lambda+vk,d+uk,\hh+\pp+uk\\e+uk,\hh+uk\end{matrix}\right)
=\frac{(-1)^{vk}\Gamma(e+\lambda-d)\Gamma(1+d-e-\lambda)\Gamma(e+uk)Y_p(u,v;k)}{(\hh+uk)_{\pp}\Gamma((u-v)k+e+\lambda)
\Gamma(vk+d-e-\lambda+p+1)},
\end{equation}
where
\begin{equation}\label{eq:Ypolynomial}
Y_p(u,v;t)=\frac{(\hh-d)_{\pp}}{\Gamma(e-d)}\sum_{j=0}^{p}\frac{(d-e+1)_{j}}{j!}
F\!\left(\begin{matrix}-j,1-\hh+d\\1-\hh+d-\pp\end{matrix}\right)(ut+d)_{j}(vt+d-e-\lambda+j+1)_{p-j}
\end{equation}
is a polynomial in $t$ of degree $p$.
\end{Lemma}
\begin{proof} According to \eqref{eq:Gpowerint} we have
\begin{multline}\label{eq:intsum}
F\!\left(\begin{matrix}-\lambda+vk,d+uk,\hh+\pp+uk\\e+uk,\hh+uk\end{matrix}\right)
\\
=\frac{\Gamma(e+uk)}{(\hh+uk)_\pp \Gamma(d+uk)}\int\limits_0^1(1-x)^{\lambda-vk}G^{p,0}_{p,p}\!\left(\!x\left|\!\begin{array}{l}e+uk-1,\hh+uk-1\\d+uk-1,\hh+\pp+uk-1\end{array}\right.\!\!\!\right)dx.
\end{multline}
Here, by \cite[Lemma~1]{KPChapter2019},
$$
G^{p,0}_{p,p}\!\left(\!x\left|\!\begin{array}{l}e+uk-1,\hh+uk-1\\d+uk-1,\hh+\pp+uk-1\end{array}\right.\!\!\!\right)
=\sum\limits_{j=0}^p\lambda_jx^{j+d+uk-1}(1-x)^{-j+e-d-1},
$$
where
$$
\lambda_j=\frac{(\hh-d)_\pp(e-d-j)_j}{\Gamma(e-d)j!}F\!\left(\begin{matrix}-j,1-\hh+d\\1-\hh+d-\pp\end{matrix}\right).
$$
Substituting this expansion into \eqref{eq:intsum} and integrating term-wise leads to
$$
F\!\left(\begin{matrix}-\lambda+vk,d+uk,\hh+\pp+uk\\e+uk,\hh+uk\end{matrix}\right)
=A\sum\limits_{j=0}^p\lambda_j\Gamma(j+d+uk)\Gamma(-j+e-d+\lambda-vk),
$$
where
$$
A=\frac{\Gamma(e+uk)}{(\hh+uk)_\pp\Gamma(d+uk)\Gamma(uk-vk+\lambda+e)}.
$$
To complete the proof, it remains to use the identities
$$
\Gamma(-j+e-d+\lambda-vk)=(-1)^{j+vk}\frac{\Gamma(e+\lambda-d)\Gamma(1-e-\lambda+d)}{\Gamma(j+1+vk-e+d-\lambda)},
$$
$$
\frac{\Gamma(j+d+uk)}{\Gamma(j+1+vk-e+d-\lambda)}=\frac{\Gamma(d+uk)}{\Gamma(1-\lambda-e+d+p+vk)}(d+uk)_j(1-\lambda-e+d+j+vk)_{p-j}.
$$
\end{proof}

\noindent\textbf{Remark.} Note that for $e=d+1$ the polynomial $Y_p(u,v;t)$ reduces to
\begin{equation}\label{eq:polKarl}
Y_p(u,v;t)=(\hh-d)_{\pp}(vt-\lambda)_{p}
\end{equation}
and \eqref{eq:IPDsummatiom} reduces to the  Karlsson-Minton summation theorem \cite[(1.3)]{KPITSF2018}
\begin{equation}\label{eq:KarlssonMinton}
F\left.\!\!\left(\!\begin{matrix}a,d,\hh+\pp\\d+1,\hh\end{matrix}\right.\right)
=\frac{\Gamma(d+1)\Gamma(1-a)(\hh-d)_{\pp}}{\Gamma(d+1-a)(\hh)_{\pp}}
\end{equation}
valid for $\Re(a+p)<1$ (recall that $p=p_1+\cdots+p_l$).

\noindent\textbf{Remark.} If $p=l=1$ the polynomial $Y_p(u,v;t)$ reduces to
\begin{equation}\label{eq:Y1}
Y_1(u,v;t)=\frac{v(h-d)-u(e-d-1)}{\Gamma(e-d)}t-\frac{(h-d)\lambda+(e-d-1)h}{\Gamma(e-d)}
\end{equation}
with the root
\begin{equation}\label{eq:Y1root}
\xi=\frac{(h-d)\lambda+(e-d-1)h}{(h-d)v-(e-d-1)u}.
\end{equation}
The cases below refer to the values of $(u,v)$ in \eqref{eq:generaltrans} and \eqref{eq:master1}.

\textbf{Case I: $(u,v)=(1,0)$.} The Pfaff-Saalsch\"{u}tz theorem \cite[Theorem~2.2.6]{AAR} in the form
$$
{}_{3}F_{2}\!\left(\begin{matrix}-n,a,b\\c,d\end{matrix}\right)
=\frac{(c-a)_{n}(d-a)_{n}}{(c)_{n}(d)_{n}},~~~~c+d=-n+a+b+1.
$$
yields
\begin{equation}\label{eq:Saalschutz}
{}_{3}F_{2}\!\left(\begin{matrix}-\lambda,-n+k,a_2+k\\b_1+k,b_2+k\end{matrix}\right)
=\frac{(b_1+\lambda)_{n}(b_2+\lambda)_{n}(b_1)_{k}(b_2)_{k}}{(b_1)_{n}(b_2)_{n}(b_1+\lambda)_{k}(b_2+\lambda)_{k}},
\end{equation}
where $b_1+b_2+n-a_2+\lambda=1$, $n\in\N$, $0\le{k}\le{n}$.

Rakha and Rathie \cite[(2.5)]{RakhaRathie} (see also \cite[(3.1)]{KRP2013})
extended the Pfaff-Saalsch\"{u}tz summation theorem by adding a parameter pair with unit shift. Their extension can be written in the form:
$$
{}_{4}F_{3}\!\left(\begin{matrix}-j,a,b,f+1\\c,d,f\end{matrix}\right)
=\frac{(c-a-1)_{j}(d-a)_{j}(\gamma+1)_{j}}{(c)_{j}(d)_{j}(\gamma)_{j}},~~~~c+d=-j+a+b+2,
$$
where
$$
\gamma=\frac{(c-a-1)(c-b-1)f}{ab+(c-a-b-1)f}.
$$
Setting $j=n-k$, $c=b_1+k$, $d=b_2+k$, $a=-\lambda$, $b=a_2+k$, $s=b_3+k$ after some rearrangements we get for $n,k\in\N$, $k\le{n}$, the following summation formula:
\begin{subequations}\label{eq:Saalunitshift}
\begin{equation}
{}_{4}F_{3}\!\left(\begin{matrix}-\lambda,-n+k,a_2+k,b_3+1+k\\b_1+k,b_2+k,b_3+k\end{matrix}\right)
=\Omega\frac{(b_1)_{k}(b_2)_{k}(b_3)_k(\mu+1)_k}
{(b_1+\lambda)_{k}(b_2+\lambda)_{k}(b_3+1)_k(\mu)_{k}},
\end{equation}
where $b_1+b_2+n-a_2+\lambda=2$,
\begin{equation}
\Omega=\frac{(na_2(b_3+\lambda)+b_3(b_1+\lambda-1)(b_2+\lambda-1))(b_1+\lambda)_{n-1}(b_2+\lambda)_{n}}{b_3(1+a_2-b_1)(b_1)_{n}(b_2)_{n}},
\end{equation}
and $\mu$ is defined by
\begin{equation}\label{eq:mu1}
\mu=\frac{na_2(b_3+\lambda)+b_3(b_1+\lambda-1)(b_2+\lambda-1)}{\lambda(b_3-b_1+1)-(b_1+n-1)(b_1-a_2-1)}.
\end{equation}
\end{subequations}
Another extension of Pfaff-Saalsch\"{u}tz's theorem is achieved by replacing $1$-balanced (or Saalsch\"{u}tzian) series with $r$-balanced series, where $r\in\N$.  The simplest formula of this type with $r=2$ as given by Kim and Rathie in \cite[(3.1)]{KR2012} can be cast into the form
$$
{}_{3}F_{2}\!\left(\begin{matrix}-j,a,b\\c,d\end{matrix}\right)
=\frac{(c-a-1)_{j}(d-a-1)_{j}}{(c)_{j}(d)_{j}}\left(1+\frac{jb}{(c-a-1)(d-a-1)}\right),~~~~c+d=-j+a+b+2.
$$
Setting $j=n-k$, $c=b_1+k$, $d=b_2+k$, $a=-\lambda$, $b=a_2+k$ after some rearrangements we get for $n,k\in\N$, $k\le{n}$, the following summation formula:
\begin{subequations}\label{eq:Saalcontig}
\begin{equation}
{}_{3}F_{2}\!\left(\begin{matrix}-\lambda,-n+k,a_2+k\\b_1+k,b_2+k\end{matrix}\right)
=B\frac{(b_1)_{k}(b_2)_{k}(\nu+1)_k}
{(b_1+\lambda)_{k}(b_2+\lambda)_{k}(\nu)_{k}},
\end{equation}
where $b_1+b_2+n+\lambda-a_2=2$
\begin{equation}
B=\frac{(b_1+\lambda)_{n-1}(b_2+\lambda)_{n-1}}{(b_1)_{n}(b_2)_{n}}(na_2+(b_1+\lambda-1)(b_2+\lambda-1)),
\end{equation}
and $\nu$ is defined by
\begin{equation}
\nu=\frac{na_2+(b_1+\lambda-1)(b_2+\lambda-1)}{b_1+b_2+n-a_2+2(\lambda-1)}.
\end{equation}
\end{subequations}

A particular case of $2$-balanced summation theorem is the following formula due to Bailey \cite[4.5(1.2)]{Bailey} (see also \cite[Table~6.1-30]{Koepf})
$$
{}_{3}F_{2}\!\left(\begin{matrix}a,b,-j\\1+a-b,1+2b-j\end{matrix}\right)
=\frac{(a-2b)_{j}(1+a/2-b)_{j}(-b)_{j}}{(1+a-b)_{j}(a/2-b)_{j}(-2b)_{j}}.
$$
Setting $b=-\lambda$, $a=\alpha+k$, $j=n-k$, we obtain:
\begin{multline}\label{eq:Bailey3F2-1}
{}_{3}F_{2}\!\left(\begin{matrix}-\lambda,\alpha+k,-n+k\\1+\lambda+\alpha+k,1-2\lambda-n+k\end{matrix}\right)
\\
=-\frac{(\alpha+2\lambda)_{n}(\lambda)_{n}(1-2\lambda-n)_{k}(1+\lambda+\alpha)_{k}(-\alpha-2\lambda-2n)_{k+1}}
{(2\lambda)_{n}(1+\lambda+\alpha)_{n}(1-\lambda-n)_{k}(-\alpha-2\lambda-2n)_{k}(\alpha+2\lambda)_{k+1}}.
\end{multline}

The $(u,v)=(1,0)$ case of Lemma~\ref{lm:IPDsummation} takes the form
\begin{equation}\label{eq:IPD10}
F\!\left(\begin{matrix}-\lambda,d+k,\hh+\pp+k\\e+k,\hh+k\end{matrix}\right)
=\frac{\Gamma(e+\lambda-d)\Gamma(e)(e)_{k}(\hh)_{k}Y_p(1,0;k)}{\Gamma(e+\lambda)(1+d-e-\lambda)_p(\hh)_{\pp}(\hh+\pp)_{k}(e+\lambda)_{k}}.
\end{equation}

\textbf{Case II: $(u,v)=(1,1)$.} In this case the only summation formula is the one given by Lemma~\ref{lm:IPDsummation}:
\begin{equation}\label{eq:IPD11}
F\!\left(\begin{matrix}-\lambda+k,d+k,\hh+\pp+k\\e+k,\hh+k\end{matrix}\right)
=\frac{(-1)^{k}\Gamma(e+\lambda-d)\Gamma(e)(e)_{k}(\hh)_{k}Y_p(1,1;k)}
{\Gamma(e+\lambda)(\hh)_{\pp}(\hh+\pp)_{k}(1+d-e-\lambda)_{p}(1+d-e-\lambda+p)_{k}}.
\end{equation}

\textbf{Case III: $(u,v)=(1,-1)$.}  Whipple's formula \cite[Table~6.1-16]{Koepf}  leads to:
\begin{align}\label{eq:Whipple16}
{}_{3}F_{2}\!\left(\begin{matrix}-\lambda-k,1+\lambda+k,a_2+k\\b_1+k,1+2a_2-b_1+k\end{matrix}\right)
=\frac{B(b_1)_{k}(1+2a_2-b_1)_{k}}
{4^k((1+\lambda+b_1)/2)_{k}((2+\lambda+2a_2-b_1)/2)_{k}},
\\[6pt]\nonumber
B=\frac{\pi{2^{1-2a_2}}\Gamma(b_1)\Gamma(1+2a_2-b_1)}
{\Gamma((b_1-\lambda)/2)\Gamma((1-\lambda+2a_2-b_1)/2)\Gamma((1+\lambda+b_1)/2)\Gamma((2+\lambda+2a_2-b_1)/2)}.
\end{align}

The $(u,v)=(1,-1)$ case of Lemma~\ref{lm:IPDsummation} reads:
$$
F\!\left(\begin{matrix}-\lambda-k,d+k,\hh+\pp+k\\e+k,\hh+k\end{matrix}\right)
=\frac{(-1)^{k}\Gamma(e+\lambda-d)\Gamma(1+d-e-\lambda)\Gamma(e+k)Y_p(1,-1;k)}{(\hh+k)_{\pp}\Gamma(2k+e+\lambda)
\Gamma(-k+d-e-\lambda+p+1)}
$$
or
\begin{equation}\label{eq:IPD1-1}
F\!\left(\begin{matrix}-\lambda-k,d+k,\hh+\pp+k\\e+k,\hh+k\end{matrix}\right)
=\frac{\Gamma(e+\lambda-d)\Gamma(e)(e)_{k}(e+\lambda-d-p)_{k}(\hh)_{k}Y_p(1,-1;k)}{\Gamma(e+\lambda)(1-e-\lambda+d)_{p}(\hh)_{\pp}(\hh+\pp)_{k}
4^k\Delta(e+\lambda,2)_{k}},
\end{equation}
where $\Delta(a,2)_{k}=(a/2)_k((a+1)/2)_k$.

\textbf{Case IV: $(u,v)=(1,2)$.}
Bailey's formula \cite[Table~6.1-30]{Koepf} for nearly-poised (of the second kind) ${}_{3}F_{2}$ is
$$
{}_{3}F_{2}\!\left(\begin{matrix}a,1+a/2,-j\\a/2,c\end{matrix}\right)
=\frac{(c-a-1-j)(c-a)_{j-1}}{(c)_{j}}.
$$
Setting $a=-\lambda+2k$, $j=n-k$, $c=b+k$, in view of $(z)_{2k}=4^k(z/2)_{k}((z+1)/2)_{k}$,
we obtain:
\begin{equation}\label{eq:Bailey2}
{}_{3}F_{2}\!\left(\begin{matrix}-\lambda+2k,1-\lambda/2+k,-n+k\\-\lambda/2+k,b+k\end{matrix}\right)
=\frac{(b+\lambda-n-1)(b+\lambda)_{n-1}(b)_{k}(1-b-\lambda)_k}{(b)_{n}(-4)^k((2-b-\lambda-n)/2)_{k}((3-b-\lambda-n)/2)_{k}}.
\end{equation}

From Dougall's formula \cite[4.3(3)]{Bailey} (see also \cite[Table~6.1-25]{Koepf})
$$
{}_{5}F_{4}\left.\!\!\left(\!\begin{matrix}a,1+a/2,c,d,-j\\a/2,1+a-c,1+a-d,1+a+j\end{matrix}\right.\right)
=\frac{(1+a)_{j}(1+a-c-d)_{j}}{(1+a-c)_{j}(1+a-d)_{j}}
$$
on setting $a=-\lambda+2k$, $c=1-\lambda-b_1+k$, $d=1-\lambda-b_2+k$, $j=n-k$ and applying the relations
$$
(\gamma)_{n-k}=\frac{(-1)^k(\gamma)_n}{(1-\gamma-n)_k},~~(\gamma+k)_{n-k}=\frac{(\gamma)_n}{(\gamma)_k},~~(z)_{2k}=4^k(z/2)_{k}((z+1)/2)_{k},
$$
we arrive at
\begin{multline}\label{eq:Dougallk}
{}_{5}F_{4}\left.\!\!\left(\!\begin{matrix}-\lambda+2k,1-\lambda/2+k,1-\lambda-b_{1}+k,1-\lambda-b_{2}+k,-n+k
\\
-\lambda/2+k,b_1+k,b_2+k,1-\lambda+n+k\end{matrix}\right.\right)
\\
=\frac{(1-\lambda)_{n}(b_1+b_2+\lambda-1)_{n}(b_1)_{k}(b_2)_{k}(1-\lambda+n)_{k}}
{(b_1)_{n}(b_2)_{n}(-4)^k(2-b_1-b_2-\lambda-n)_{k}((1-\lambda)/2)_{k}((2-\lambda)/2)_{k}}.
\end{multline}

The $(u,v)=(1,2)$ case of Lemma~\ref{lm:IPDsummation} by application of
\begin{equation}\label{eq:GammaTrans}
\Gamma(2z)=2^{2z-1}\pi^{-1/2}\Gamma(z)\Gamma(z+1/2),~~~ \Gamma(z-k)=\frac{(-1)^k\Gamma(z)}{(1-z)_{k}}
\end{equation}
takes the form
\begin{multline}\label{eq:IPD12}
F\!\left(\begin{matrix}-\lambda+2k,d+k,\hh+\pp+k\\e+k,\hh+k\end{matrix}\right)
\\
=\frac{\Gamma(e+\lambda-d)\Gamma(e)(\hh)_{k}(e)_{k}(1-e-\lambda)_{k}Y_p(1,2;k)}
{\Gamma(e+\lambda)(1+d-e-\lambda)_{p}(\hh)_{\pp}(\hh+\pp)_{k}(-4)^k\Delta(1-e-\lambda+p+d,2)_{k}}.
\end{multline}

\textbf{Case V: $(u,v)=(2,2)$.} The $(u,v)=(2,2)$ case of Lemma~\ref{lm:IPDsummation} after application of \eqref{eq:GammaTrans} reads:
\begin{multline}\label{eq:IPD22}
F\!\left(\begin{matrix}-\lambda+2k,d+2k,\hh+\pp+2k\\e+2k,\hh+2k\end{matrix}\right)
\\
=\frac{\Gamma(e+\lambda-d)\Gamma(e)\Delta(e,2)_{k}\Delta(\hh,2)_{k}Y_p(2,2;k)}
{\Gamma(e+\lambda)(1-e-\lambda+d)_{p}(\hh)_{\pp}\Delta(\hh+\pp,2)_{k}\Delta(1-e-\lambda+p+d,2)_{k}}.
\end{multline}

\textbf{Case VI: $(u,v)=(2,1)$.} The $(u,v)=(2,1)$ case of Lemma~\ref{lm:IPDsummation} after application of \eqref{eq:GammaTrans} reads:
\begin{equation}\label{eq:IPD21}
F\!\left(\begin{matrix}-\lambda+k,d+2k,\hh+\pp+2k\\e+2k,\hh+2k\end{matrix}\right)
=\frac{\Gamma(e+\lambda-d)\Gamma(e)(-4)^k\Delta(e,2)_{k}\Delta(\hh,2)_{k}Y_p(2,1;k)}
{(\hh)_{\pp}\Gamma(e+\lambda)\Delta(\hh+\pp,2)_{k}(e+\lambda)_{k}(1+d-e-\lambda)_{p+k}}.
\end{equation}
Another formula that can sum the hypergeometric function on the right hand side of (8) in Case VI can be obtained from Watson's formula \cite[Table 6.1-16]{Koepf}. However, its application does not lead to any new or interesting known transformations, so we omit it here.

\subsection{Transformation formulas}
In this subsection we present a number of known transformation formulas of the type given in \eqref{eq:generaltrans}.  
We will group them into the same six cases that we have used in the previous subsection.
The cases $(u,v)=(1,1)$ and $(u,v)=(1,2)$ have subcases with $w=1$ and $w=2$, but as the value of $w$ does not affect the summation on the right hand side of \eqref{eq:master1} we keep these two situations under the single case.

\textbf{Case I: $(u,v)=(1,0)$}. The second Euler-Pfaff transformation  belongs to this class:
\begin{equation}\label{eq:Euler-2}
{}_{2}F_{1}\left.\!\!\left(\begin{matrix}a, b
\\c\end{matrix}\right\vert x\right)
=(1-x)^{c-a-b}{}_{2}F_{1}\left.\!\!\left(\begin{matrix}c-a,
c-b\\c\end{matrix}\right\vert x\right).
\end{equation}
Its natural extension to the hypergeometric functions with integral  parameter differences is the second  Miller-Paris transformation \cite{MP2013,KPChapter2019}.
Define  $\m=(m_1,\ldots,m_r)\in\N^r$, $m=m_1+m_2+\ldots+m_r$ and  $\f=(f_1,\ldots,f_r)\in\C^{r}$.  We will reserve he symbols $\f$ and $\m$ for the Miller-Paris transformations and their corollaries throughout the rest of the paper.
According to \cite[Theorem~1]{KRP2014} and \cite[Theorem~4]{MP2013} we have:
\begin{equation}\label{eq:KRPTh1-2}
F\left.\!\left(\begin{matrix}a, b,\f+\m\\c,\f\end{matrix}\right\vert x\right)
=(1-x)^{c-a-b-m}F\left.\!\left(\begin{matrix}c-a-m,c-b-m, \hat{\zetta}+1\\c, \hat{\zetta}\end{matrix}\right\vert x\right),
\end{equation}
where $\hat{\zetta}=(\zeta_1,\ldots,\zeta_m)$ are the roots of the characteristic polynomial
\begin{equation}\label{eq:Qmhat}
\hat{Q}_m(a,b,c,\f,\m;t)=\sum\limits_{k=0}^{m}\frac{{}_{r+1}F_{r}\!\left(-k,\f+\m;\f\right)(a)_k(b)_k(t)_k}{(c-a-m)_k(c-b-m)_kk!}
{}_{3}F_{2}\!\left(\begin{matrix}-m+k,t+k,c-a-b-m\\c-a-m+k,c-b-m+k\end{matrix}\right).
\end{equation}
Transformation \eqref{eq:KRPTh1-2} holds for $|x|<1$ when $(1+a+b-c)_m\ne0$, $(c-a-m)_m\ne0$ and  $(c-b-m)_m\ne0$.
A somewhat simpler but less symmetric form of the characteristic polynomial $\hat{Q}_m$ was given by us in \cite[(15)]{KPChapter2019}.

\textbf{Case II: $(u,v)=(1,1)$}.  The first Euler-Pfaff transformation is given by
\begin{equation}\label{eq:Euler-1}
{}_{2}F_{1}\left.\!\!\left(\begin{matrix}a, b
\\c\end{matrix}\right\vert x\right)
=(1-x)^{-a}{}_{2}F_{1}\left.\!\!\left(\begin{matrix}a, c-b
\\c\end{matrix}\right\vert\frac{x}{x-1}\right).
\end{equation}
It was extended to hypergeometric functions with integral  parameter differences
by Miller and Paris \cite{MP2013,KPChapter2019}. Keeping the meaning of $\m$ and $\f$ we have according to \cite[Theorem~1]{KRP2014} and \cite[Theorem~3]{MP2013}:
\begin{equation}\label{eq:KRPTh1-1}
F\left.\!\left(\begin{matrix}a, b,\f+\m\\c,\f\end{matrix}\right\vert x\right)
=(1-x)^{-a}F\left.\!\left(\begin{matrix}a,c-b-m, \zetta+1\\c,\zetta\end{matrix}\right\vert\frac{x}{x-1}\right)
\end{equation}
where $\zetta=\zetta(c,b,\f)=(\zeta_1,\ldots,\zeta_m)$ are the roots of the polynomial
\begin{equation}\label{eq:Qm}
Q_m(b,c,\f,\m;t)=\frac{1}{(c-b-m)_{m}}\sum\limits_{k=0}^{m}\frac{(-1)^k}{k!}{}_{r+1}F_{r}\!\left(\begin{matrix}-k,\f+\m\\\f\end{matrix}\right)(b)_k(t)_{k}(c-b-m-t)_{m-k}.
\end{equation}
Two alternative forms of this polynomial can be found in \cite[(3.7)]{KPITSF2018} and \cite[Theorem~1]{KPChapter2019}.  Transformation \eqref{eq:KRPTh1-1}
is valid when $b\neq f_j,$ $j=1,\ldots,r$ and  $(c-b-m)_m\neq 0$. Both formulas  \eqref{eq:KRPTh1-2} and \eqref{eq:KRPTh1-1} fail
when $c-b-m\in\{-m+1,\ldots,0\}$. We called this situation degenerate and found the extensions of \eqref{eq:KRPTh1-2} and
\eqref{eq:KRPTh1-1} to this case in our recent papers \cite{KPChapter2019,KPResults2019}.
Miller-Paris transformations reduce to Euler-Pfaff transformations  \eqref{eq:Euler-2}, \eqref{eq:Euler-1} when $m=0$.

Further, according to \cite[(3.3)]{KrRao2003} for $|x|<1$:
\begin{equation}\label{eq:quadKrRao3.3}
{}_{2}F_{1}\!\!\left(\left.\!\begin{matrix}\alpha,\alpha+1/2\\\beta\end{matrix}\right|x^2\right)
=(1-x)^{-2\alpha}{}_{2}F_{1}\!\!\left(\left.\!\begin{matrix}2\alpha,\beta-1/2
\\2\beta-1\end{matrix}\right|\frac{-2x}{1-x}\right).
\end{equation}

\textbf{Case III: $(u,v)=(1,-1)$}.  According to \cite[(3.5)]{KrRao2003} for $x<1/2$ we have the Gauss transformation
\begin{equation}\label{eq:quadKrRao3.5}
{}_{2}F_{1}\!\!\left(\left.\!\begin{matrix}\alpha,\beta\\(\alpha+\beta+1)/2\end{matrix}\right|x\right)
={}_{2}F_{1}\!\!\left(\left.\!\begin{matrix}\alpha/2,\beta/2
\\(\alpha+\beta+1)/2\end{matrix}\right|4x(1-x)\right).
\end{equation}
This formula remains true for all $x$ if both sides terminate ($\alpha$ and/or $\beta$ is a negative even integer).

Further, according to \cite[(3.9)]{KrRao2003} we have for $x<1/2$
\begin{equation}\label{eq:quadKrRao3.9}
{}_{2}F_{1}\!\!\left(\left.\!\begin{matrix}\alpha,1-\alpha\\\beta\end{matrix}\right|x\right)
=(1-x)^{\beta-1}{}_{2}F_{1}\!\!\left(\left.\!\begin{matrix}(\beta-\alpha)/2,(\alpha+\beta-1)/2
\\\beta\end{matrix}\right|4x(1-x)\right).
\end{equation}
This formula remains true for all $x$ if both sides terminate (for instance, when $\alpha>\beta$ are positive integers both odd or both even).

\textbf{Case IV: $(u,v)=(1,2)$}.
The following transformation is known as Kummer's first quadratic transformation \cite[(6.6)]{MP2013} (cf. \cite[15.8.15]{NIST}):
\begin{equation}\label{eq:quad15.8.15}
{}_{2}F_{1}\!\!\left(\left.\!\begin{matrix}\alpha,\beta\\1-\beta+\alpha\end{matrix}\right|-x\right)
=(1-x)^{-\alpha}{}_{2}F_{1}\!\!\left(\left.\!\begin{matrix}\alpha/2,\alpha/2+1/2
\\1-\beta+\alpha\end{matrix}\right|\frac{-4x}{(1-x)^2}\right).
\end{equation}
It is true for $|x|<1$.

Next, according to \cite[(3.1.11)]{AAR} we have for $|x|<1$
\begin{equation}\label{eq:quad15.8.21}
{}_{2}F_{1}\!\!\left(\left.\!\begin{matrix}\alpha,\beta\\1-\beta+\alpha\end{matrix}\right|x^2\right)
=(1-x)^{-2\alpha}{}_{2}F_{1}\!\!\left(\left.\!\begin{matrix}\alpha,\alpha-\beta+1/2
\\1-2\beta+2\alpha\end{matrix}\right|\frac{-4x}{(1-x)^2}\right).
\end{equation}

Whipple's quadratic transformation \cite[(3.1.15)]{AAR} is given by
\begin{equation}\label{eq:quadWhipple}
{}_{3}F_{2}\!\!\left(\left.\!\begin{matrix}\alpha,\beta,\delta\\1-\beta+\alpha,1-\delta+\alpha\end{matrix}\right|x\right)
=(1-x)^{-\alpha}{}_{3}F_{2}\!\!\left(\left.\!\begin{matrix}\alpha/2,(\alpha+1)/2,1+\alpha-\beta-\delta
\\1-\beta+\alpha, 1-\delta+\alpha\end{matrix}\right|\frac{-4x}{(1-x)^2}\right),
\end{equation}
which is also valid  for $|x|<1$.

According to Choi and Rathie \cite[(2.1)]{ChoiRathie} (after change of variable and change of notation), we have
\begin{equation}\label{eq:quadChRat2.1}
{}_{2}F_{1}\!\!\left(\left.\!\begin{matrix}\alpha,\beta\\\beta+1\end{matrix}\right|x\right)
=(1-x)^{-2\beta}{}_{3}F_{2}\!\!\left(\left.\!\begin{matrix}\beta,\beta-\alpha/2+1,\beta-\alpha/2+1/2
\\\beta+1, 2\beta-\alpha+1\end{matrix}\right|\frac{-4x}{(1-x)^2}\right).
\end{equation}
We will refer to the above transformation as the first Choi-Rathie transformation.
A closely related result \cite[(2.4)]{ChoiRathie} after change of variable and change of notation takes the form:
\begin{equation}\label{eq:quadChRat2.4}
{}_{3}F_{2}\!\!\left(\left.\!\begin{matrix}\alpha+1,2\alpha,\beta\\\alpha,\beta+1\end{matrix}\right|x\right)
=(1-x)^{-2\beta}{}_{3}F_{2}\!\!\left(\left.\!\begin{matrix}\beta,\beta-\alpha,\beta-\alpha+1/2
\\\beta+1, 2\beta-2\alpha+1\end{matrix}\right|\frac{-4x}{(1-x)^2}\right),
\end{equation}
where  $|x|<1$. We will refer to this  transformation as the second Choi-Rathie transformation.

A result by Rakha and Rathie \cite[(3.1)]{RakhaRathie} reads
\begin{equation}\label{eq:quadRakhaRathie3.1}
{}_{4}F_{3}\!\left(\left.\begin{matrix}2\alpha,\alpha-\beta-1/2,1+\alpha-\sigma,1+\alpha+\sigma
\\\alpha+\beta+3/2,\alpha-\sigma,\alpha+\sigma\end{matrix}\right|x\right)
=(1-x)^{-2\alpha}{}_{3}F_{2}\!\left(\left.\!\begin{matrix}\alpha,\beta,\delta+1
\\\alpha+\beta+3/2,\delta\end{matrix}\right|\frac{-4x}{(1-x)^2}\right)
\end{equation}
for $|x|<1$, where
$$
\sigma^2=\frac{1}{\beta-\delta}\left(\alpha^2\beta-\alpha\beta\delta-\beta\delta/2-\delta/4\right)
$$
with $\sigma^2<0$ permitted.

The following more recent transformation is given by Wang and Rathie in \cite[(3.1)]{WangRathie}:
\begin{multline}\label{eq:quadWangRathie3.1}
{}_{5}F_{4}\!\left(\left.\!\begin{matrix}2\alpha-1,2\alpha-\beta-1,2\alpha-\gamma,1/2+\alpha-\omega,1/2+\alpha+\omega
\\\beta+1,\gamma,\alpha-1/2-\omega,\alpha-1/2+\omega,\end{matrix}\:\right|x\right)
\\
=(1-x)^{1-2\alpha}{}_{4}F_{3}\!\left(\left.\!\begin{matrix}\alpha,\alpha-1/2,\beta+\gamma-2\alpha,\delta+1
\\\beta+1,\gamma,\delta\end{matrix}\right|\frac{-4x}{(1-x)^2}\right)
\end{multline}
for $|x|<1$, where
$$
\omega^2=\left(\alpha-1/2\right)^2-\frac{\delta(\gamma-2\alpha)(2\alpha-\beta-1)}{\beta+\gamma-2\alpha-\delta}.
$$

Kummer's first transformation \eqref{eq:quad15.8.15} was generalized by Miller
and Paris in \cite[Theorem~5]{MP2013} to the generalized hypergeometric functions with integral parameter differences as follows:
\begin{equation}\label{eq:MP2013-6.3}
F\!\left(\left.\!\begin{matrix}\alpha,\beta,\hat{\etta}+1\\1-\beta+\alpha,\hat{\etta}\end{matrix}\:\right|-x\right)
=(1-x)^{-\alpha}F\!\left(\left.\!\begin{matrix}\alpha/2,\alpha/2+1/2,\f+\m
\\1-\beta+\alpha,\f\end{matrix}\right|\frac{-4x}{(1-x)^2}\right)
\end{equation}
where $|x|<1$ and $\hat{\etta}$ is the vector of zeros of the $2m$ degree polynomial
\begin{equation}\label{eq:hatR2m}
\hat{R}_{2m}(t;\alpha,\beta,\f,\m)=\sum_{k=0}^{m}\frac{(\f)_{\m}(t)_{k}(\alpha-t)_{k}}{(\beta)_kk!}{}_{r+1}F_{r}\!\!\left(\!\begin{matrix}-k,\f+\m\\\f\end{matrix}\right).
\end{equation}

Another set of extensions of the classical quadratic transformations has been obtained recently by Maier \cite{Maier2019}.
Whipple's transformation \eqref{eq:quadWhipple} is generalized to \cite[Theorem~3.1]{Maier2019}
\begin{equation}\label{eq:MaierTh3.1}
F\!\left(\left.\!\begin{matrix}\alpha,\beta,\delta,1-\rrho\\1+\alpha-\beta,1+\alpha-\delta,-\rrho\end{matrix}\:\right|x\right)
=(1-x)^{-\alpha}{}_{3}F_{2}\!\left(\left.\!\begin{matrix}\alpha/2,\alpha/2+1/2,\alpha-\beta-\delta-k+1
\\
1+\alpha-\beta,1+\alpha-\delta\end{matrix}\right|\frac{-4x}{(1-x)^2}\right),
\end{equation}
where $\rrho$ is the vector of roots of the $2k$ degree polynomial
\begin{equation}\label{eq:MaierQ2kpoly}
P_{2k}(t;\alpha,\beta,\delta)={}_{3}F_{2}\!\left(\!\begin{matrix}-k,-t,t+\alpha\\\beta,\delta\end{matrix}\right).
\end{equation}
Further extension has been obtained by adding a parameter pair $\left[\begin{matrix}\gamma+k\\\gamma\end{matrix}\right]$ on the right hand side \cite[Theorem~3.4]{Maier2019}:
\begin{multline}\label{eq:MaierTh3.4}
F\!\left(\left.\!\begin{matrix}\alpha,\beta,\delta,1-\hat{\rrho}\\1+\alpha-\beta,1+\alpha-\delta,-\hat{\rrho}\end{matrix}\:\right|x\right)
\\
=(1-x)^{-\alpha}{}_{4}F_{3}\!\left(\left.\!\begin{matrix}\alpha/2,\alpha/2+1/2,\alpha-\beta-\delta-k+1,\gamma+k
\\
1+\alpha-\beta,1+\alpha-\delta,\gamma\end{matrix}\right|\frac{-4x}{(1-x)^2}\right),
\end{multline}
where $\hat{\rrho}$ is the vector of roots of the $2k$ degree  polynomial
\begin{equation}\label{eq:MaierhatQ2kpoly}
\hat{P}_{2k}(t;\alpha,\beta,\delta,\gamma)={}_{4}F_{3}\!\left(\!\begin{matrix}-k,-t,t+\alpha,\beta+\delta+\gamma+k-\alpha-1\\\beta,\delta,\gamma\end{matrix}\right).
\end{equation}
Renaming parameters it is easy to see that \eqref{eq:MaierTh3.4} is a generalization of \eqref{eq:quadWangRathie3.1} to which it reduces when $k=1$.
One more extension is given in \cite[Theorem~3.7]{Maier2019}, but with characteristic polynomial defined recursively. We omit this case here.

\textbf{Case V: $(u,v)=(2,2)$}.
The following transformation is known as Kummer's second quadratic transformation \cite[(6.5)]{MP2013} (cf. \cite[15.8.13]{NIST}):
\begin{equation}\label{eq:quad15.8.13}
{}_{2}F_{1}\!\!\left(\left.\!\begin{matrix}\alpha,\beta\\2\beta\end{matrix}\:\right|2x\right)
=(1-x)^{-\alpha}{}_{2}F_{1}\!\!\left(\left.\!\begin{matrix}\alpha/2,\alpha/2+1/2
\\\beta+1/2\end{matrix}\right|\frac{x^2}{(1-x)^2}\right),
\end{equation}
where For $-1<x<1/2$.  This formula remains valid for $-1<x<1$ if we assume that $-\alpha\in\N$ so that both sides terminate.

Kummer's transformation \eqref{eq:quad15.8.13} was generalized by Miller and Paris in \cite[Theorem~5]{MP2013} to
generalized hypergeometric functions with integral parameter differences as follows. According to \cite[(6.1)]{MP2013} we have:
\begin{equation}\label{eq:MP2013-6.1}
F\!\left(\left.\!\begin{matrix}\alpha,\beta-m,\etta+1\\2\beta,\etta\end{matrix}\right|2x\right)
=(1-x)^{-\alpha}F\!\!\left(\left.\!\begin{matrix}\alpha/2,\alpha/2+1/2,\f+\m
\\\beta+1/2,\f\end{matrix}\right|\frac{x^2}{(1-x)^2}\right)
\end{equation}
for $x<1/2$, where $\etta$ is the vector of roots of the polynomial
\begin{equation}\label{eq:R2m}
R_{2m}(t)=\sum_{k=0}^{m}\frac{(-1)^k(\f)_{\m}(t)_{2k}(\beta-m-t)_{m-k}}{4^kk!}{}_{r+1}F_{r}\!\!\left(\!\begin{matrix}-k,\f+\m\\\f\end{matrix}\right).
\end{equation}
Formula \eqref{eq:MP2013-6.1} is true for $0<x<1$ if we assume that $-\alpha\in\N$ so that both sides terminate.

\textbf{Case VI: $(u,v)=(2,1)$}. According to \cite[15.8.14]{NIST} we have for $|x|<1$
\begin{equation}\label{eq:quad15.8.14}
{}_{2}F_{1}\!\!\left(\left.\!\begin{matrix}\alpha,\beta\\2\beta\end{matrix}\right|x\right)
=(1-x)^{-\alpha/2}{}_{2}F_{1}\!\!\left(\left.\!\begin{matrix}\alpha/2,\beta-\alpha/2
\\\beta+1/2\end{matrix}\right|\frac{-x^2}{4(1-x)}\right).
\end{equation}

\section{$G$ function integral method: results}

By application of Lemma 1 to case $i$ transformation, $i\in\{I, II, III, IV, V, VI\}$, playing the role of \eqref{eq:generaltrans}, and using case $i$ summation formulas for summing the generalized hypergeometric function on the right hand side of \eqref{eq:master1} we arrive at the transformation formulas below grouped according the values of $(u,v)$ in \eqref{eq:generaltrans}.  

A remark is here in order regarding the convergence regions of the identities presented below. According to \cite[Theorem~2.1.2]{AAR} the generalized hypergeometric series \eqref{eq:pFqdefined} with $P(n)\equiv1$ converges absolutely at $x=1$ if $\Re(\sum{b_i}-\sum{a_j})>0$.
If $P(n)$ is a polynomial of degree $m$, it following immediately from the definition of the Pochhammer symbol or also from \eqref{eq:FP-roots}   that this condition must be modified to $\Re(\sum{b_i}-\sum{a_j})>m$. This condition gives the convergence regions for the identities involving non-terminating series. For terminating series we have finite summations so that the identities are true for all values of parameters such that no denominator vanishes. 

\subsection{Case I: $(u,v)=(1,0)$}
Fix $r,l\in\N$, $\m\in\N^r$, $\f\in\C^{r}$, $\pp\in\N^l$, $\hh\in\C^{l}$, $a,b,d,e\in\C$ (recall that $m=m_1+\cdots+m_r$, $p=p_1+\cdots+p_{l}$).  By an application of  the beta integral method to the second Miller-Paris
transformation \eqref{eq:KRPTh1-2}, Kim, Rathie and Paris proved in \cite{KRP2014} that
\begin{multline}\label{eq:KRPar1}
F\left.\!\!\left(\!\begin{matrix}a,b,d,\f+\m,\hh+\pp\\c,e,\f,\hh\end{matrix}\right.\right)
\\
=\frac{\Gamma(e)\Gamma(c+e-a-b-d-m-p)}{\Gamma(e-d)\Gamma(c+e-a-b-m-p)}
{}F\left.\!\!\left(\!\begin{matrix}c-a-m-p,c-b-m-p,d,\hat{\zetta}_*+1\\c,c+e-a-b-m-p,\hat{\zetta}_*\end{matrix}\right.\right),
\end{multline}
where $\Re(e-d)>0$, $\Re(c+e-a-b-d-m-p)>0$, $(c-b-m-p)_{m+p}\ne0$, $(c-b-a)_{m+p}\ne0$,
$\hat{\zetta}_*$ are the roots of the polynomial $\hat{Q}_{m+p}(a,b,c,(\f,\hh),(\m,\pp);t)$ defined in \eqref{eq:Qmhat}.

The following  theorem shows that the above formula can be viewed as an extreme case of a family of
transformations of the left hand side.
\begin{Theorem}\label{th:MPIPD}
Suppose that $(1+a+b-c)_m\ne0$, $(c-a-m)_m\ne0$, $(c-b-m)_m\ne0$ and convergence conditions $\Re(e-d-p)>0$ and $\Re(c+e-a-b-d-m-p)>0$ are satisfied. Then
\begin{equation}\label{eq:masterIPD}
F\left.\!\!\left(\!\begin{matrix}a,b,d,\f+\m,\hh+\pp\\c,e,\f,\hh\end{matrix}\right.\right)
=\Lambda{\cdot}F\left.\!\left(\!\begin{array}{l}c-a-m,c-b-m,d,\hat{\zetta}+1 \\c,e+c-a-b-m,\hat{\zetta} \end{array}\, \vdots\, Y_p(1,0)\right.\right),
\end{equation}
where
$$
\Lambda=\frac{\Gamma(e+c-a-b-m-d)\Gamma(e)}{{(\hh)_\pp}\Gamma(e+c-a-b-m)(1+d-e-c+a+b+m)_p},
$$
the polynomial $Y_p(1,0)=Y_p(1,0;z)$ is defined in \eqref{eq:Ypolynomial} with $\lambda=c-a-b-m$, $\hat{\zetta}$ are the roots of the polynomial $\hat{Q}_{m}(a,b,c,\f,\m;t)$ defined in \eqref{eq:Qmhat}.

Formula \eqref{eq:masterIPD} remains valid for $m=0$, if the parameters $\f+\m$, $\f$, $\hat{\zetta}+1$, $\hat{\zetta}$ are omitted.
\end{Theorem}

\begin{proof}  Conditions of the theorem ensure that transformation \eqref{eq:KRPTh1-2} holds. This transformation is a particular case of \eqref{eq:generaltrans} if we identify the parameters as follows:
$$
\aalpha=(a,b,\f+\m),~~~\bbeta=(c,\f),~~\lambda=c-a-b-m,~~\ddelta=(c-a-m,c-b-m,\hat{\zetta}+1),
$$
$$
\ggamma=(c,\hat{\zetta}),~~D=1,~~~u=1, ~~~v=0.
$$
Setting $\a=(d,\hh+\pp)$, $\b=(e,\hh)$ we can apply  Lemma~\ref{lm:master1} to conclude that
\begin{multline*}
F\left.\!\!\left(\!\begin{matrix}a,b,d,\f+\m,\hh+\pp\\c,e,\f,\hh\end{matrix}\right.\right)
\\
=\sum\limits_{k=0}^{\infty}\frac{(c-a-m)_k(c-b-m)_k(\hat{\zetta}+1)_k(d)_{k}(\hh+\pp)_k}{(c)_k(e)_{k}(\hh)_k(\hat{\zetta})_k k!}
F\left.\!\!\left(\!\begin{matrix} a+b+m-c,d+k,\hh+\pp+k\\e+k,\hh+k\end{matrix}\right.\right).
\end{multline*}
Summing the hypergeometric function of the right hand side by formula \eqref{eq:IPD10} we arrive at \eqref{eq:masterIPD}.
If $m=0$, instead of the Miller-Paris transformation  \eqref{eq:KRPTh1-2} start with
the second Euler-Pfaff transformation \eqref{eq:Euler-2}. 
\end{proof}

If $e=d+1$ the polynomial $Y_p(1,0;z)$ has the form \eqref{eq:polKarl}, so that  \eqref{eq:masterIPD} reduces to yet another extension of Karlsson-Minton summation theorem \eqref{eq:KarlssonMinton}:
$$
F\left.\!\!\left(\!\begin{matrix}a,b,d,\f+\m,\hh+\pp\\c,d+1,\f,\hh\end{matrix}\right.\!\right)
\!=\!
\frac{\Gamma(c-a-b-m+1)\Gamma(d+1)(\hh-d)_\pp}{{}\Gamma(d+c-a-b-m+1)(\hh)_\pp}
F\left.\!\!\left(\!\!\!\begin{array}{l}c-a-m,c-b-m,d,\hat{\zetta}+1 \\c,d+c-a-b-m+1,\hat{\zetta} \end{array}\right.\!\!\!\right).
$$
This formula holds provided that $c-a-m$ or $c-b-m$ is a negative integer.

For $r=l=1,$ $\f=(f),$ $\hh=(h),$ $m_1=p_1=1$ the polynomial
$\hat{Q}_1(t)$ takes the form \cite[p.116]{KRP2014}
$$
\hat{Q}_1(t)=1+\frac{(c-a-b-1)f+ab}{(c-a-1)(c-b-1)}\frac{t}{f}.
$$
Then, in view of \eqref{eq:Y1}, \eqref{eq:Y1root}, formula \eqref{eq:masterIPD} reduces to
\begin{multline}\label{eq:5F4-2unitshifts}
{}_{5}F_{4}\left.\!\!\left(\!\begin{matrix}
a,b,d,h+1,f+1\\c,e,h,f\end{matrix}\right.\right)
=
\frac{((e-d-1)h+(c-a-b-1)(h-d))\Gamma(e)\Gamma(s_*)}{h\Gamma(s_*+d+1)\Gamma(e-d)}
\\
\times{}_{5}F_{4}\left.\!\!\left(\!\begin{matrix}
c-a-1,c-b-1,d,\hat{\xi}+1,\hat{\zeta}+1\\c,e+c-a-b-1,\hat{\xi},\hat{\zeta}\end{matrix}\right.\right),
\end{multline}
where $s_*=e+c-a-b-d-2$, $c-a-1\ne0$, $c-b-1\ne0$,
$$
\hat{\xi}=h+\frac{(c-a-b-1)(h-d)}{e-d-1},~~~~\hat{\zeta}=\frac{(c-a-1)(c-b-1)f}{(c-a-b-1)f+ab}.
$$
Note that $\hat{\xi}$ and $\hat{\zeta}$ are linear-fractional functions of parameters, while, in contrast,  application of the Kim, Rathie, Paris formula \eqref{eq:KRPar1} to the left hand side of \eqref{eq:5F4-2unitshifts} leads to ${}_5F_{4}$ on the right hand side containing the conjugate quadratic roots among parameters.  Setting $d=h$ leads to $m=1$ case of \eqref{eq:KRPar1}.

Setting $m=0$, $l=1$,  $\hh=(h)$, $\pp=(1)$. Then formula
\eqref{eq:masterIPD} from Theorem~\ref{th:MPIPD} takes the form
\begin{multline}\label{eq:4F3generate}
{}_{4}F_{3}\left.\!\!\left(\!\begin{matrix}
a,b,d,h+1\\c,e,h\end{matrix}\right.\right)=
\frac{((e-d-1)h+(c-a-b)(h-d))\Gamma(e)\Gamma(s)}{(e-d-1)h\Gamma(s+d+1)\Gamma(e-d-1)}
{}_{4}F_{3}\left.\!\!\left(\!\begin{matrix}
c-a,c-b,d,\xi_*+1
\\
c,e+c-a-b,\xi_*\end{matrix}\right.\right),
\end{multline}
where  $s=e+c-a-b-d-1$, $\xi_*=h+(c-a-b)(h-d)/(e-d-1)$. We remark that formula \eqref{eq:KRP4F3unitshift} obtained by setting $m=0$, $p=1$ in \eqref{eq:KRPar1} (see \cite[p.116]{KRP2014}) has the right hand side essentially different from the one above.  Both our identity above and \eqref{eq:KRP4F3unitshift} can be applied to themselves repeatedly.  We found several other transformations connecting the ${}_4F_{3}$ functions with one unit shift and undertook a group-theoretic study of their properties in \cite{KPMath2020}. The group-theoretic properties
of terminating Saalsch\"{u}tzian  ${}_{4}F_{3}$ (i.e. with parametric excess equal to unity) have been studied
in \cite{KrRao2004,RaoBook,RDN}.

\begin{Theorem}
Suppose that $(1+a+b-c)_m\ne0$, $(c-a-m)_m\ne0$, $(c-b-m)_m\ne0$,  $n\in\N$ and
$g+e+n-d+\lambda=1$, where $\lambda=c-a-b-m$. Then
\begin{equation}\label{eq:p2Saalschutz}
F\!\left(\begin{matrix}-n,a,b,d,\f+\m\\c,g,e,\f\end{matrix}\right)
=\frac{(g+\lambda)_n(e+\lambda)_n}{(g)_n(e)_n}F\!\left(\begin{matrix}-n,a+\lambda,b+\lambda,d,\hat{\zetta}+1
\\c,g+\lambda,e+\lambda,\hat{\zetta}\end{matrix}\right),
\end{equation}
where $\hat{\zetta}$ is the vector of zeros of the polynomial
$\hat{Q}_m(a,b,c,\f,\m;t)$ defined in
\eqref{eq:Qmhat}.  Hypergeometric functions on both sides of the above formula are Saalsch\"{u}tzian.
\end{Theorem}

\begin{proof}  Put $\a=(-n,d)$, $\b=(g,e)$. Apply Lemma \ref{lm:master1} to the  transformation \eqref{eq:KRPTh1-2} playing the role of \eqref{eq:generaltrans} (with parameter identification $\aalpha=(a,b,\f+\m)$, $\bbeta=(c,\f)$, $\lambda=c-a-b-m$, $\ddelta=(c-a-m,c-b-m,\hat{\zetta}+1)$, $\ggamma=(c,\hat{\zetta})$, $D=1$, $u=1,$ $v=0$).  Use formula  (\ref{eq:Saalschutz}) to sum the hypergeometric function on the right hand side to complete the proof. \end{proof}

The most useful case of the above theorem is $r=m=1$:
\begin{multline}\label{eq:5F4Saal}
{}_{5}F_{4}\!\left(\begin{matrix}-n,a,b,d,f+1\\c,g,e,f\end{matrix}\right)
\\
=\frac{(g+c-a-b-1)_n(e+c-a-b-1)_n}{(g)_n(e)_n}{}_{5}F_{4}\!\left(\begin{matrix}-n,c-a-1,c-b-1,d,\hat{\zeta}+1\\c,g+c-a-b-1,e+c-a-b-1,\hat{\zeta}\end{matrix}\right),
\end{multline}
where $g+e+n-d+c-a-b=2$ and
$$
\hat{\zeta}=\frac{(c-a-1)(c-b-1)f}{(c-a-b-1)f+ab}.
$$

Letting $n\to\infty$, $g=2-n+d-e-c+a+b\to-\infty$, while keeping other parameters fixed,  we obtain
\begin{equation}\label{eq:KRP4F3unitshift}
{}_{4}F_{3}\!\left(\begin{matrix}a,b,d,f+1\\c,e,f\end{matrix}\right)
=\frac{\Gamma(e)\Gamma(e+c-a-b-1-d)}{\Gamma(e-d)\Gamma(e+c-a-b-1)}{}_{4}F_{3}\!\left(\begin{matrix}c-a-1,c-b-1,d,\hat{\zeta}+1\\c,e+c-a-b-1,\hat{\zeta}\end{matrix}\right)
\end{equation}
- a particular case of the Kim, Rathie and Paris formula \eqref{eq:KRPar1} derived by the beta integral method \cite[p.116]{KRP2014}.  The limit transition can be justified by Tannery's theorem which is a particular case of the Lebesgue dominated convergence theorem.

If we let $a,c\to\infty$ while $c-a=2-n+d+b-g-e$ is fixed, we arrive at a transformation for general terminating ${}_4F_{3}$ with one unit shift:
$$
{}_{4}F_{3}\!\left(\begin{matrix}-n,d,b,f+1\\g,e,f\end{matrix}\right)
=\frac{(g-d)_n(e-d)_n}{(g)_n(e)_n}{}_{4}F_{3}\!\left(\begin{matrix}-n,d,1-n+d+b-g-e,\hat{\zeta}_*+1\\1-g+d-n,1-e+d-n,\hat{\zeta}_*\end{matrix}\right),
$$
where $\hat{\zeta}_*=f(1+b+d-n-g-e)/b$.

If we let $f\to\infty$ in \eqref{eq:5F4Saal} we obtain (recall that $\lambda=c-a-b-1$):
$$
{}_{4}F_{3}\!\left(\begin{matrix}-n,d,a,b\\g,e,c\end{matrix}\right)
=\frac{(g+\lambda)_n(e+\lambda)_n}{(g)_n(e)_n}{}_{5}F_{4}\!\left(\begin{matrix}-n,d,a+\lambda,b+\lambda,\hat{\zeta}_{**}+1\\g+\lambda,e+\lambda,c,\hat{\zeta}_{**}\end{matrix}\right),
$$
where
$\hat{\zeta}_{**}=(c-a-1)(c-b-1)/(c-a-b-1)$
and the condition $g+e+c+n-a-b-d=2$ must be satisfied. This condition says that the ${}_{4}F_{3}$ on the left hand side is $2$-balanced, while ${}_5F_4$ on the right hand side is Saalsch\"{u}tzian.  As the right hand side above can be written as a linear combination of two ${}_4F_3$ functions, this formula can be viewed as a
three-term relation for terminating  $2$-balanced ${}_4F_3$.

Setting $d=f$, $g=f+1$ in \eqref{eq:5F4Saal} we get a
Saalsch\"{u}tzian ${}_3F_{2}$ on the left hand side.  The condition
$g+e+c+n-a-b-d=2$ becomes $e+c-a-b-1=-n$, so that the function
${}_{5}F_{4}$ on the right hand side reduces to  ${}_{4}F_{3}$ truncated at the $n$-term.  Using the notation $[{}_{4}F_{3}]_{n}$
for such truncated series and renaming the parameters according to
$A=f$, $B=c-a-1,$ $D=c-b-1$, $E=c$, $G=f+c-a-b$,
 we get the following curious summation formula
\begin{equation}\label{eq:4F3truntedsum}
\left[{}_{4}F_{3}\!\left(\begin{matrix}A,B,D,\hat{\zeta}+1\\G,E,\hat{\zeta}\end{matrix}\right)\right]_{n}
=\frac{(A+1)_n(B+1)_n(D+1)_n}{(G)_n(E)_nn!},
\end{equation}
where the formula for $\hat{\zeta}$ in terms of the new parameters takes the
form:
$$
\hat{\zeta}=\frac{e_3(A,B,D)}{e_2(A,B,D)-e_2(1-G,1-E)}.
$$
Here $e_k$ is $k$-th elementary symmetric polynomial. The
Saalsch\"{u}tzian condition
$e_{1}(A,B,D)+e_{1}(1-G,1-E)=0$ must be satisfied for
the validity of this formula. This is a summation formula with \emph{non-linearly constrained parameters} - a rather rare species in the hypergeometric literature.  Letting $n\to\infty$ in this formula
we recover our recent result \cite[(45)]{KPResults2019}.  For instance, if $A=B=C$, $3A=-2(1-G)=-2(1-E)$ we get
$$
\sum\limits_{j=0}^{n}\frac{[(A)_j]^3(A+3j/4)}{[(3A/2+1)_j]^2j!}=\frac{A[(A+1)_n]^3}{[(3A/2+1)_n]^2n!}.
$$

On the other hand, if we set $d=f$, $c=f+1$ in
\eqref{eq:5F4Saal} we get Saalsch\"{u}tzian ${}_3F_{2}$ on both
sides which does not lead to any new formulas.

\begin{Theorem}\label{th:SuperSaalschutz}
Suppose that $(1+a+b-c)_m\ne0$, $(c-a-m)_m\ne0$, $(c-b-m)_m\ne0$,  $n\in\N$ and
$g+e+n-d+\lambda=2$, where $\lambda=c-a-b-m$. Then
\begin{equation}\label{eq:SuperSaalschutz}
F\!\left(\begin{matrix}-n,a,b,d,h+1,\f+\m\\c,g,e,h,\f\end{matrix}\right)
=\Omega{\cdot}F\!\left(\begin{matrix}-n,d,a+\lambda,b+\lambda,\mu+1,\hat{\zetta}+1\\g+\lambda,e+\lambda,c,\mu,\hat{\zetta}\end{matrix}\right),
\end{equation}
where
\begin{equation}\label{eq:Omega}
\Omega=\frac{(g+\lambda)_{n-1}(e+\lambda)_{n-1}}{h(g)_{n}(e)_{n}}(nd(h+\lambda)+h(g+\lambda-1)(e+\lambda-1)),
\end{equation}
\begin{equation}\label{eq:mu2}
\mu=\frac{nd(h+\lambda)+h(g+\lambda-1)(e+\lambda-1)}{\lambda(h-g+1)-(g+n-1)(g-d-1)}
\end{equation}
and $\hat{\zetta}$ is the vector of zeros of the polynomial
$\hat{Q}_m(a,b,c,\f,\m;t)$ defined in \eqref{eq:Qmhat}.  Hypergeometric functions on both sides of are \eqref{eq:SuperSaalschutz} Saalsch\"{u}tzian.
\end{Theorem}

\begin{proof}   Set $\a=(-n,d,h+1)$,  $\b=(g,e,h)$.   We again apply  Lemma~\ref{lm:master1} to the  transformation \eqref{eq:KRPTh1-2} (with parameters identified as follows:
$\aalpha=(a,b,\f+\m),$ $\bbeta=(c,\f)$, $\lambda=c-a-b-m,$ $\ddelta=(c-a-m,c-b-m,\hat{\zetta}+1)$, $\ggamma=(c,\hat{\zetta})$, $D=1$, $u=1,$ $ v=0$).
The right hand  side of  \eqref{eq:master1} from Lemma~\ref{lm:master1} is then
$$
\sum\limits_{k=0}^{\infty}\frac{(\ddelta)_k(-n)_k(d)_k(h+1)_k}{(\ggamma)_k(e)_k(g)_k(h)_kk!}
{}_{4}F_{3}\left.\!\left(\!\begin{matrix}-\lambda,-n+k,d+k,h+k+1\\g+k,e+k,h+k\end{matrix}\right.\right).
$$
It remains to apply \eqref{eq:Saalunitshift} to sum the hypergeometric function on the right hand side. 
\end{proof}

The most useful case of the above Theorem is $r=m=1$:
\begin{equation}\label{eq:6F5Saal}
{}_{6}F_{5}\!\left(\begin{matrix}-n,d,a,b,h+1,f+1\\g,e,c,h,f\end{matrix}\right)
=\Omega\cdot{}_{6}F_{5}\!\left(\begin{matrix}-n,d,a+\lambda,b+\lambda,\mu+1,\hat{\zeta}+1\\g+\lambda,e+\lambda,c,\mu,\hat{\zeta}\end{matrix}\right),
\end{equation}
where $\hat{\zeta}=(c-a-1)(c-b-1)f/[ab+(c-a-b-1)f]$, $\lambda=c-a-b-1$ and  Saalsch\"{u}tz's condition $g+e+n-d+\lambda=2$ is satisfied. The numbers $\Omega$ and $\mu$ are given in \eqref{eq:Omega} and \eqref{eq:mu2}, respectively.

If we let $n\to\infty$, $g=2-n+d-e-\lambda\to-\infty$, while keeping other parameters fixed,  we recover formula
\eqref{eq:5F4-2unitshifts}.

Next, if $f\to\infty$ in \eqref{eq:6F5Saal} we obtain (recall that $\lambda=c-a-b-1$):
$$
{}_{5}F_{4}\!\left(\begin{matrix}-n,d,a,b,h+1\\g,e,c,h\end{matrix}\right)
=\Omega\cdot{}_{6}F_{5}\!\left(\begin{matrix}-n,d,a+\lambda,b+\lambda,\mu+1,\hat{\zeta}_*+1\\g+\lambda,e+\lambda,c,\mu,\hat{\zeta}_*\end{matrix}\right),
$$
where $\hat{\zeta}_*=(c-a-1)(c-b-1)/(c-a-b-1)$, $\Omega$ and $\mu$ are defined by \eqref{eq:Omega} and \eqref{eq:mu2},
respectively, and the condition $g+e+n-d+\lambda=2$  is satisfied.  This condition states that the ${}_6F_{5}$ on the right hand side is Saalsch\"{u}tzian, while  while ${}_{5}F_{4}$ on the left hand side is $2$-balanced.

Finally, if we let $a,c\to\infty$ while $c-a=3-n+d+b-g-e$ is fixed, we arrive at the transformation
$$
{}_{5}F_{4}\!\left(\begin{matrix}-n,b,d,h+1,f+1\\g,e,h,f\end{matrix}\right)
=\hat{\Omega}
\cdot{}_{5}F_{4}\!\left(\begin{matrix}-n,d,2-g-e-n+d+b,\hat{\mu}+1,\hat{\zeta}_{**}+1\\2-g-n+d,2-e-n+d,\hat{\mu},\hat{\zeta}_{**}\end{matrix}\right),
$$
where $\hat{\zeta}_{**}=f(2+b+d-n-g-e)/b$,
$$
\hat{\mu}=\frac{nd(2+h+d-n-e-g)+h(1+d-n-e)(1+d-n-g)}{(2+d-n-e-g)(h-g+1)-(g+n-1)(g-d-1)}
$$
and
\begin{multline*}
\hat{\Omega}=\frac{(2+d-n-e)_{n-1}(2+d-n-g)_{n-1}}{h(g)_{n}(e)_{n}}
\\
\times(nd(2+h+d-n-e-g)+h(1+d-n-e)(1+d-n-g)).
\end{multline*}
The function on the right hand side has the same type as the function on the left hand side (terminating ${}_{5}F_{4}$ with two unit shifts), so that this transformation can be iterated.

The proofs of the above theorems follow the same simple pattern: an application of Lemma~\ref{lm:master1} to a Case I transformation followed by an application of a suitable summation formula. Therefore, below we simply list the remaining results obtained in this way for the case $(u,v)=(1,0)$.

Combining \eqref{eq:KRPTh1-2} with \eqref{eq:Saalcontig} (renaming $a_2=d$, $b_1=g,$ $b_2=e$) we obtain the transformation
$$
F\left.\!\!\left(\!\begin{matrix}-n,a,b,d,\f+\m\\c,g,e,\f\end{matrix}\right.\right)
=B{\cdot}F\left.\!\!\left(\!\begin{matrix}-n,a+\lambda,b+\lambda,d,\nu+1,\hat{\zetta}+1\\c,g+\lambda,e+\lambda,\nu,\hat{\zetta}\end{matrix}\right.\right),
$$
where  $\lambda=c-a-b-m,$ $g+e+n+\lambda-d=2$, $(1+a+b-c)_m\ne0$, $(c-a-m)_m\ne0$,  $(c-b-m)_m\ne0$,  $\hat{\zetta}$ is the vector of zeros of the polynomial $\hat{Q}_m(a,b,c,\f,\m;t)$ and according to \eqref{eq:Saalcontig}
$$
B=\frac{(g+\lambda)_{n-1}(e+\lambda)_{n-1}}{(g)_{n}(e)_{n}}(nd+(g+\lambda-1)(e+\lambda-1)),
$$
%\begin{equation}\label{eq:nu}
$$
\nu=\frac{(nd+(g+\lambda-1)(e+\lambda-1))}{g+e+n-d+2(\lambda-1)}.
$$
%\end{equation}

Combining \eqref{eq:KRPTh1-2} with \eqref{eq:Bailey3F2-1} we obtain the transformation
%\begin{equation}\label{eq:master1-11}
$$
F\left.\!\!\left(\!\begin{matrix}-n,a,b,\alpha,\f+\m\\c,1+\lambda+\alpha,1-2\lambda-n,\f\end{matrix}\right.\right)
=G{\cdot}F\left.\!\!\left(\!\begin{matrix}-n,a+\lambda,b+\lambda,\alpha,-\alpha-2\lambda-2n+1,\hat{\zetta}+1\\c,1-\lambda-n,\alpha+2\lambda+1,-\alpha-2\lambda-2n,\hat{\zetta}\end{matrix}\right.\right),
$$
%\end{equation}
where $\lambda=c-a-b-m$, $(1+a+b-c)_m\ne0$, $(c-a-m)_m\ne0$,  $(c-b-m)_m\ne0$,  $\hat{\zetta}$ is the vector of zeros of the polynomial $\hat{Q}_m(a,b,c,\f,\m;t)$ defined in \eqref{eq:Qmhat} and
$$
G=-\frac{(\alpha+2\lambda)_{n}(\lambda)_{n}(-\alpha-2\lambda-2n)}
{(2\lambda)_{n}(1+\lambda+\alpha)_{n}(\alpha+2\lambda)}.
$$

\subsection{Case II: $(u,v)=(1,1)$}

Fix $r,l\in\N$, $\m\in\N^r$, $\f\in\C^{r}$, $\pp\in\N^l$, $\hh\in\C^{l}$, $a,b,d,e\in\C$ (recall that $m=m_1+\cdots+m_r$, $p=p_1+\cdots+p_{l}$).  By an application of  the beta integral method to the first Miller-Paris
transformation \eqref{eq:KRPTh1-1}, Kim, Rathie and Paris proved in \cite{KRP2014} that
\begin{equation}\label{eq:KRPar2}
F\left.\!\!\left(\!\begin{matrix}-n,b,d,\f+\m,\hh+\pp\\c,e,\f,\hh\end{matrix}\right.\right)
=\frac{(e-d)_n}{(e)_n}
F\left.\!\!\left(\!\begin{matrix}-n,c-b-m-p,d,\zetta_*+1\\c,1-e+d-n,\zetta_*\end{matrix}\right.\right),
\end{equation}
where  $\zetta_*$ are the roots of the polynomial
${Q}_{m+p}(b,c,(\f,\hh),(\m,\pp);t)$, $\Re(e-d)>0$, $b\neq f_j$, $1\le{j}\le{r}$, $(c-b-m-p)_{m+p}\ne0$.
Similarly to the previous case $(u,v)=(1,0)$, our approach embeds this identity into a family of transformations.  Members of this family generally have two characteristic polynomials of lower degree: one of degree $m$ and the other of degree $p$ in contrast to one polynomial ${Q}_{m+p}$ of degree $m+p$ for \eqref{eq:KRPar2}.

\begin{Theorem}\label{thm:2}
Suppose $b\ne f_j$, $j=1,\ldots,r$, $(c-b-m)_m\ne0$. Then for each $n\in\N$ we have
\begin{equation}\label{eq:masterIPD-2}
F\left.\!\!\left(\!\begin{matrix}-n,b,d,\f+\m,\hh+\pp\\c,e,\f,\hh\end{matrix}\right.\right)
=
\frac{(e-d)_n\Gamma(e-d)}{(e)_n(1+d-e-n)_p(\hh)_\pp} F\left.\!\left(\!\begin{array}{l}-n,c-b-m,d,\zetta+1 \\c,1+d-e-n+p,\zetta \end{array}\, \vdots \, Y_p(1,1)\, \right.\right).
\end{equation}
Here  the polynomial $Y_p(1,1)=Y_p(1,1;z)$ is defined in \eqref{eq:Ypolynomial} with $\lambda=n$, $\zetta$ is the vector of zeros of the polynomial ${Q}_{m}(a,b,c,\f,\m;t)$ defined in \eqref{eq:Qm}.

Formula \eqref{eq:masterIPD-2} remains valid in the case $m = 0$, if we omit the parameters $\f+\m,$ $\f,$ $\zetta+1,$ $\zetta.$
\end{Theorem}

\begin{proof}  Conditions of the theorem imply that transformation \eqref{eq:KRPTh1-1} holds true.  This transformation is a particular case of \eqref{eq:generaltrans} if we identify the parameters as follows: $w=1$, $M=1$, $\aalpha=(-n,b,\f+\m)$, $\bbeta=(c,\f)$,  $\lambda=n$, $\ddelta=(-n,c-b-m,\zetta+1)$, $\ggamma=(c,\zetta)$, $D=-1$, $u=1$, $v=1$.  Setting  $\a=(d,\hh+\pp)$, $\b=(e,\hh)$ we apply  Lemma~\ref{lm:master1} to get
\begin{multline*}
F\left.\!\!\left(\!\begin{matrix}-n,b,d,\f+\m,\hh+\pp\\c,e,\f,\hh\end{matrix}\right.\right)
\\
=\sum\limits_{k=0}^{\infty}\frac{(-n)_k(c-b-m)_k(\zetta+1)_k(d)_{k}(\hh+\pp)_k(-1)^k}{(c)_k(e)_{k}(\hh)_k(\zetta)_k k!}
F\left.\!\!\left(\!\begin{matrix} -n+k,d+k,\hh+\pp+k\\e+k,\hh+k\end{matrix}\right.\right).
\end{multline*}
Application of the summation formula \eqref{eq:IPD11} to the hypergeometric function on the right hand side completes the proof.  In the case $m=0$, instead of the first Miller-Paris transformation  \eqref{eq:KRPTh1-1} start with the first  Euler-Pfaff transformation \eqref{eq:Euler-1}.
\end{proof}

For $r=l=1,$ $\f=(f),$ $\hh=(h),$ $m_1=p_1=1$, $\lambda=n$ by using \eqref{eq:FP-roots} formula \eqref{eq:masterIPD-2} from Theorem~\ref{thm:2} takes the form
$$
{}_{5}F_{4}\left.\!\!\left(\!\begin{matrix}
-n,b,d,h+1,f+1\\c,e,h,f\end{matrix}\right.\right)
=\frac{(e-d)_n}{h(e)_n}\left(h+\frac{nd}{1+d-e-n}\right)
{}_{5}F_{4}\left.\!\!\left(\!\begin{matrix}
-n,c-b-1,d,\zeta+1,\xi_*+1\\c,2+d-e-n,\zeta,\xi_*\end{matrix}\right.\right),
$$
where
$$
\xi_*=\frac{(h-d)n+(e-d-1)h}{e-h-1},~~~\zeta=\frac{(c-b-1)f}{f-b}.
$$
Here $\xi_*$ is the negated root of $Y_1(1,1;z)$ according to \eqref{eq:Y1root} and $\zeta$ is is the root of  $Q_1(t)$ according to \eqref{eq:Qm} (see also \cite[p.116]{KRP2014}). We further applied
\eqref{eq:Y1} to express $Y_1(1,1;0)$.

In a similar fashion, setting $m=0$, $l=1$,  $\hh=(h)$, $\pp=(1)$,  formula \eqref{eq:masterIPD-2} takes the form
\begin{equation*}\label{eq:th2-4}
{}_{4}F_{3}\left.\!\!\left(\!\begin{matrix}-n,b,d,h+1
\\
c,e,h\end{matrix}\right.\right)
=\frac{(e-d)_n}{(e)_n}\left(1+\frac{dn}{h(1+d-e-n)}\right)
{}_{4}F_{3}\left.\!\!\left(\!\begin{matrix} -n,c-b,d,\xi_*+1
\\
c, 2+d-e-n,\xi_*\end{matrix}\right.\right),
\end{equation*}
where $\xi_*$ is as defined above.  We note that formula \eqref{eq:KRPar2} due to Kim, Rathie, Paris
under these restrictions reads \cite[p.116]{KRP2014}
\begin{equation*}\label{eq:th1-5}
{}_{4}F_{3}\left.\!\!\left(\!\begin{matrix}
-n,b,d,f+1\\c,e,f\end{matrix}\right.\right)=
\frac{(e-d)_n}{(e)_n}
{}_{4}F_{3}\left.\!\!\left(\!\begin{matrix}
-n,c-b-1,d,\zeta+1\\c,1-e-d-n,\zeta\end{matrix}\right.\right),
\end{equation*}
where $\zeta=(c-b-1)f/(f-b)$.

Finally, combination \eqref{eq:quadKrRao3.3} with \eqref{eq:IPD11} yields a rather unusual  transformation involving a terminating hypergeometric function evaluated at $2$:
\begin{equation}\label{eq:master1-12}
F\left.\!\!\left(\!\begin{matrix}\alpha,\alpha+1/2,\Delta(d,2),\Delta(\hh+\pp,2)\\\beta,\Delta(e,2),\Delta(\hh,2)\end{matrix}\right.\right)
=C\cdot F\left.\!\!\left(\!\begin{array}{l}2\alpha,\beta-1/2,d,e\\
2\beta-1,e,1+d-e+2\alpha+p\end{array}\, \vdots \, Y_p(1,1)\, \right|2\right),
\end{equation}
where $-d\in\N$,
$$
C=\frac{\Gamma(e-2\alpha-d)\Gamma(e)}
{\Gamma(e-2\alpha)(\hh)_{\pp}(1+d-e+2\alpha)_{p}},
$$
 and the polynomial $Y_p(1,1)$ is defined in \eqref{eq:Ypolynomial} with $\lambda=-2\alpha$.

\subsection{Case III: $(u,v)=(1,-1)$}

Combining \eqref{eq:quadKrRao3.5} with \eqref{eq:Whipple16} and renaming parameters according to $A=\alpha$, $B=\beta$, $C=a_2$, $D=b_1$  we obtain the transformation
$$
{}_{4}F_{3}\left.\!\!\left(\!\begin{matrix}1,A,B,C
\\
(A+B+1)/2,D,1+2C-D\end{matrix}\right.\right)
={}_{4}F_{3}\left.\!\!\left(\!\begin{matrix}1,A/2,B/2,C
\\
(A+B+1)/2,(1+D)/2,1+C-D/2\end{matrix}\right.\right),
$$
valid if both sides terminate.

Combining \eqref{eq:quadKrRao3.5} with \eqref{eq:IPD1-1} and writing $a=\alpha$, $b=\beta$ we obtain the transformation
$$
F\left.\!\!\left(\!\begin{matrix}a,b,d,\hh+\pp
\\
(a+b+1)/2,e,\hh\end{matrix}\right.\right)
=\frac{\Gamma(e-d)}{(1-e+d)_{p}(\hh)_{\pp}}
F\!\left(\!\begin{array}{l}a/2,b/2,d,e-d-p\\(a+b+1)/2,\Delta(e,2)\end{array}\, \vdots \, Y_p(1,-1)\, \right),
$$
where $Y_p(1,-1)$ is defined in \eqref{eq:Ypolynomial} with $\lambda=0$, and both sides must terminate. This transformation can further be extended to any values of parameters making both sides converge using Carlson's theorem (see an example of such extension in Case IV below).  If $p=l=1$ it reduces to
$$
{}_{4}F_{3}\left.\!\!\left(\!\begin{matrix}a,b,d,h+1
\\
(a+b+1)/2,e,h\end{matrix}\right.\right)
={}_{5}F_{4}\!\left(\!\begin{array}{l}a/2,b/2,d,e-d-1,\xi_{*}+1\\(a+b+1)/2,\Delta(e,2),\xi_{*}\end{array}\!\!\right),
$$
where, according to \eqref{eq:Y1root},
$$
\xi_{*}=-\frac{(e-d-1)h}{2d-h-e+1}.
$$

Combining \eqref{eq:quadKrRao3.9} with \eqref{eq:IPD1-1} and writing $a=\alpha$, $b=\beta$  we obtain the transformation
\begin{multline}\label{KrRao3.9-IPD}
F\left.\!\!\left(\!\begin{matrix}a,1-a,d,\hh+\pp
\\
b,e,\hh\end{matrix}\right.\right)
\\
=\frac{\Gamma(e+b-d-1)\Gamma(e)}{\Gamma(e+b-1)(2-e-b+d)_{p}(\hh)_{\pp}}
F\!\left(\!\begin{array}{l}(b-a)/2,(a+b-1)/2,d,e+b-d-p-1\\b,\Delta(e+b-1,2)\end{array}\, \vdots \, Y_p(1,-1)\, \right),
\end{multline}
where $Y_p(1,-1)$ is defined in \eqref{eq:Ypolynomial} with $\lambda=b-1$, and both sides must terminate.
If $p=l=1$ it reduces to
\begin{multline*}
{}_{4}F_{3}\left.\!\!\left(\!\begin{matrix}a,1-a,d,h+1
\\
b,e,h\end{matrix}\right.\right)
=\left(1+\frac{d(b-1)}{h(2-e-b+d)}\right)
\frac{\Gamma(e+b-d-1)\Gamma(e)}{\Gamma(e+b-1)\Gamma(e-d)}
\\
\times{}_{5}F_{4}\!\left(\!\begin{array}{l}(b-a)/2,(a+b-1)/2,d,e+b-d-2,\xi_{**}+1
\\
b,(e+b-1)/2,(e+b)/2,\xi_{**}\end{array}\!\!\right),
\end{multline*}
where both sides must terminate and, according to \eqref{eq:Y1root},
$$
\xi_{**}=-\frac{(h-d)(b-1)+(e-d-1)h}{2d-h-e+1}.
$$

\subsection{Case IV: $(u,v)=(1,2)$}
The following  transformations are obtained by combining case IV transformations with case IV summation formulas. Their proofs follow the same simple pattern which we illustrate by giving a proof of the first transformation. All subsequent formulas are proved in a similar fashion.

\begin{enumerate}

%1
\item Combination of Kummer's first transformation \eqref{eq:quad15.8.15} with \eqref{eq:Dougallk}
leads to a transformation of the general very well-poised ${}_6F_{5}(-1)$ to ${}_3F_{2}(1)$:
\begin{multline}\label{eq:6F5(-1)-VWP-3F24-1Bal}
{}_{6}F_{5}\left.\!\!\left(\!\begin{matrix}A,1+A/2,B,C,D,E
\\
A/2,1+A-B,1+A-C,1+A-D,1+A-E\end{matrix}\right|-1\right)
=\frac{\Gamma(1+A-C-D-E)}{\Gamma(1+A)}
\\
\times\frac{\Gamma(1+A-C)\Gamma(1+A-D)\Gamma(1+A-E)}{\Gamma(1+A-C-D)\Gamma(1+A-D-E)\Gamma(1+A-C-E)}
{}_{3}F_{2}\left.\!\!\left(\!\begin{matrix}C,D,E
\\
1+A-B,C+D+E-A\end{matrix}\right.\right),
\end{multline}
where $-E\in\N$.  This is easily seen to be equivalent to Bailey's formula
\cite[4.4(2)]{Bailey} (see also \cite[Theorem~3.4.6]{AAR}) by an application of the Thomae's relation
\cite[section~3.2]{Bailey} to  ${}_{3}F_{2}$ to the RHS.  The restriction $-E\in\N$ is then removed by the
fact that Bailey's formula is $n\to\infty$ limit of Whipple's transformation \eqref{eq:7F6-VWP-4F3-1bal},
so that the above identity remains true if parameters are restricted to make both sides converge.

For the proof of \eqref{eq:6F5(-1)-VWP-3F24-1Bal} apply Lemma~\ref{lm:master1} to the particular case of \eqref{eq:generaltrans}  given by \eqref{eq:quad15.8.15}
to get
$$
F\!\left(\left.\!\begin{matrix}\alpha,\beta,\a\\1-\beta+\alpha,\b\end{matrix}\right|-1\right)
=\sum\limits_{k}\frac{(\alpha/2)_k(\alpha/2+1/2)_k(-4)^k(\a)_k}{(1-\beta+\alpha)_{k}k!(\b)_k}
F\left.\!\!\left(\!\begin{matrix}\alpha+2k,\a+k
\\
\b+k\end{matrix}\right.\right).
$$
Next we choose the hypergeometric function on the right hand side to fit the summation formula \eqref{eq:Dougallk} (with $\lambda=-\alpha$), i.e. setting $a_1=1+\alpha/2$, $a_2=1+\alpha-b_{1}$, $a_3=1+\alpha-b_{2}$, $a_4=-n$, $n\in\N$, $b_{0}=\alpha/2$, with arbitrary $b_{1}$, $b_{2}$  and $b_{3}=1+\alpha+n$.  Application of \eqref{eq:Dougallk} after some cancelations and renaming the parameters according to $A=\alpha$, $B=\beta$, $C=1+\alpha-b_1$, $D=1+\alpha-b_2$, $E=-n$  leads to \eqref{eq:6F5(-1)-VWP-3F24-1Bal}.

%2
\item Combination of \eqref{eq:quad15.8.21} with \eqref{eq:Dougallk} after renaming the parameters according to  $A=\alpha$, $B=\beta$, $C=1+2\alpha-b_1$, $D=1+2\alpha-b_2$, $E=-n$, yields a presumably new transformation connecting a particular case of well-poised ${}_9F_{8}$ to ${}_4F_{3}$ which is neither balanced nor well-poised:
\begin{multline}\label{eq:9F8WPspecial}
{}_{9}F_{8}\left.\!\!\left(\!\begin{matrix}A,1+A/2,B,\Delta(C,2),\Delta(D,2),\Delta(E,2)
\\
A/2,1+A-B,\Delta(1+2A-C,2),\Delta(1+2A-D,2),\Delta(1+2A-E,2)
\end{matrix}\right.\right)
\\
=\frac{\Gamma(1+2A-C)\Gamma(1+2A-D)\Gamma(1+2A-E)\Gamma(1+2A-C-D-E)}
{\Gamma(1+2A)\Gamma(1+2A-C-D)\Gamma(1+2A-C-E)\Gamma(1+2A-D-E)}
\\
\times{}_{4}F_{3}\left.\!\!\left(\!\begin{matrix}A-B+1/2,C,D,E
\\
A+1/2,2A-2B+1,C+D+E-2A\end{matrix}\right.\right),
\end{multline}
where $-E\in\N$.   We will prove that this formula remains true if the series on the left hand side converges,
while the series on the right hand side terminates.  The proof is by an application of Carlson's theorem \cite[p.40]{Bailey}.
Indeed,  writing $E=-n-z$, $n\in\N$, we have proved the above identity for $z=0,1,2,\ldots$ Next, assume that the parameters are restricted so that the ${}_7F_{6}$ series
obtained on the left hand side by deleting the  parameters containing $E$ is convergent,
i.e
$$
\Re(1+2A-B-C-D)>0.
$$
The terms containing $E$ take the form
$$
\frac{(-n/2-z/2)_{k}(-n/2-z/2+1/2)_{k}}{(1/2+A+n/2+z/2)_{k}(1+A+n/2+z/2)_{k}}
$$
which is uniformly (in $k$) bounded for $\Re(z)\ge0$ if $\Re(1+2A)>0$.  Under these restrictions
the function  on the left hand side is holomorphic and bounded in the half-plane $\Re(z)\ge0$.
The function
$$
\frac{\Gamma(1+2A+n+z)\Gamma(1+2A-C-D+n+z)}{\Gamma(1+2A-C+n+z)\Gamma(1+2A-D+n+z)}
$$
on the right hand side is holomorphic and bounded in $\Re(z)\ge0$ if we additionally assume that $\Re(1+2A-C-D)>0$.
Finally, the series ${}_4F_{3}$ on the right hand side consists of a finite number of terms, say $M$, and
has poles at the points:
$$
z=C+D-2A-n+j,~~~j=0,1,\ldots,M
$$
All these points lie in the left half-plane if $M<\Re(2A+n-C-D)$ and each terms is bounded under this
condition.  Hence, for any finite $M$ we can find sufficiently large $n$ in order that the above
condition be satisfied.  We are then in the position to apply Carlson's theorem to conclude that both
sides are equal for $\Re(z)\ge0$.  Additional assumptions made above can now be removed by
analytic continuation.

%3
\item Combination of \eqref{eq:quadWhipple} with \eqref{eq:Dougallk} gives (after renaming parameters) Whipple's transformation
\cite[Theorem~3.4.4]{AAR} of very well-poised ${}_7F_{6}$ to 1-balanced ${}_4F_{3}$:
\begin{multline}\label{eq:7F6-VWP-4F3-1bal}
{}_{7}F_{6}\left.\!\!\left(\!\begin{matrix}A,1+A/2,B,C,D,E,F
\\
A/2,1+A-B,1+A-C,1+A-D,1+A-E,1+A-F\end{matrix}\right.\right)
\\
=\frac{\Gamma(1+A-D)\Gamma(1+A-E)\Gamma(1+A-F)\Gamma(1+A-D-E-F)}{\Gamma(1+A)\Gamma(1+A-D-F)\Gamma(1+A-D-E)\Gamma(1+A-E-F)}
\\
\times{}_{4}F_{3}\left.\!\!\left(\!\begin{matrix}1+A-B-C,D,E,F
\\
1+A-B,1+A-C,D+E+F-A\end{matrix}\right.\right),
\end{multline}
which is valid for $-F\in\N$.  The formula is then extended to any values of parameters such that the left hand sides converges while the right hand side terminates using Carlson's theorem. See details in \cite[section~5.4]{Bailey}.

%4
\item Combination of \eqref{eq:quadChRat2.1} with \eqref{eq:Dougallk} gives transformation connecting general nearly poised (of the first kind) ${}_4F_{3}$ with  1-balanced ${}_5F_{4}$ discovered by Bailey \cite[p.32, 4.6(1)]{Bailey}. After renaming parameters it takes the form
\begin{equation}\label{eq:4F3-NP-5F4-1bal}
{}_{4}F_{3}\left.\!\!\left(\!\begin{matrix}A,B,C,-n
\\
\kappa-B,\kappa-C,\kappa+n\end{matrix}\right.\right)
=\frac{(\kappa)_n(\kappa-B-C)_{n}}{(\kappa-B)_n(\kappa-C)_{n}}
{}_{5}F_{4}\left.\!\!\left(\!\begin{matrix}\Delta(\kappa-A,2),B,C,-n
\\
\Delta(\kappa,2),\kappa-A,B+C-\kappa+1-n\end{matrix}\right.\right),
\end{equation}
where $n\in\N$.

%5
\item Combination of \eqref{eq:quadChRat2.4}  with \eqref{eq:Dougallk} after renaming parameters according to $A=2\alpha$, $\kappa=1+2\beta$, $B=1+2\beta-b_1$, $C=1+2\beta-b_2$, $D=-n$, leads to a transformation of a particular nearly-poised (of the first kind) ${}_5F_{4}$ to a particular 2-balanced ${}_5F_{4}$:
\begin{multline}\label{eq:5F4-NP-5F4-2bal}
{}_{5}F_{4}\left.\!\!\left(\!\begin{matrix}A,1+A/2,B,C,D
\\
A/2,\kappa-B,\kappa-C,\kappa-D\end{matrix}\right.\right)
=\frac{\Gamma(\kappa-B)\Gamma(\kappa-C)\Gamma(\kappa-D)\Gamma(\kappa-B-C-D)}
{\Gamma(\kappa)\Gamma(\kappa-B-C)\Gamma(\kappa-B-D)\Gamma(\kappa-C-D)}
\\
\times{}_{5}F_{4}\left.\!\!\left(\!\begin{matrix}\Delta(\kappa-A-1,2),B,C,D
\\
\Delta(\kappa,2),\kappa-A,B+C+D-\kappa+1\end{matrix}\right.\right),
\end{multline}
where $-D\in\N$.  This relation resembles \cite[p.32, 4.6(2)]{Bailey}, but does not reduce to it.
Parameters could be extended to cover the case when the left hand side converges while the right side terminates using Carlson's theorem. Furthermore, numerical experiments show that this identity remains true for any parameters making both sides convergent.

%6
\item Combination of \eqref{eq:quad15.8.15} and \eqref{eq:Bailey2} gives transformation of the second kind nearly-poised   ${}_4F_3(-1)$
to another ${}_4F_{3}(1)$ which neither poised nor balanced:
\begin{multline}\label{eq:4F3-NP2-4F3}
{}_{4}F_{3}\left.\!\!\left(\!\begin{matrix}A,1+A/2,B,-n
\\
A/2,1-B+A,C\end{matrix}\right|-1\right)
=\frac{(C-A-n-1)(C-A)_{n-1}}{(C)_{n}}
\\
\times{}_{4}F_{3}\left.\!\!\left(\!\begin{matrix}1+A/2,(1+A)/2,1-C+A,-n
\\
1-B+A,(2-C+A-n)/2,(3-C+A-n)/2\end{matrix}\right.\right),
\end{multline}
where $n\in\N$.

%7
\item Combination of \eqref{eq:quad15.8.21} and \eqref{eq:Bailey2} gives the following transformation:
\begin{multline}\label{eq:5F4special-4F3}
{}_{5}F_{4}\left.\!\!\left(\!\begin{matrix}A,1+A/2,B,-n/2,(-n+1)/2
\\
A/2,1-B+A,C/2,(C+1)/2\end{matrix}\right.\right)
=\frac{(C-2A-n-1)(C-2A)_{n-1}}{(C)_{n}}\times
\\
{}_{4}F_{3}\left.\!\!\left(\!\begin{matrix}1+A,A-B+1/2,1-C+2A,-n
\\
1-2B+2A,(2-C+2A-n)/2,(3-C+2A-n)/2\end{matrix}\right.\right),
\end{multline}
where $n\in\N$.

%8
\item Combination of  Whipple's quadratic transformation \eqref{eq:quadWhipple} with \eqref{eq:Bailey2} gives the following transformation  of the second kind nearly-poised  ${}_5F_{4}$ to Saalsch\"{u}tzian ${}_5F_{4}$ which was discovered by Bailey \cite[4.5(2)]{Bailey}:
\begin{multline}\label{eq:5F4-NP2-5F4-1bal}
{}_{5}F_{4}\left.\!\!\left(\!\begin{matrix}A,1+A/2,B,C,-n
\\
A/2,1-B+A,1-C+A,D\end{matrix}\right.\right)
=\frac{(D-A-n-1)(D-A)_{n-1}}{(D)_{n}}\times
\\
{}_{5}F_{4}\left.\!\!\left(\!\begin{matrix}1+A/2,(1+A)/2,1-D+A,1+A-B-C,-n
\\
1-B+A,1-C+A,(2-D+A-n)/2,(3-D+A-n)/2\end{matrix}\right.\right),
\end{multline}
where $n\in\N$.

%9
\item Combination of the Rakha-Rathie transformation \eqref{eq:quadRakhaRathie3.1}  with Dougall's summation formula \eqref{eq:Dougallk} leads to a transformation of a particular Saalsch\"{u}tzian ${}_5F_{4}$ with one unit shift to very well-poised ${}_8F_{7}$ with two unit shifts. Renaming the parameters according to $A=2\alpha$, $B=\alpha-\beta-1/2$, $C=1+2\alpha-b_1$, $D=1+2\alpha-b_2$, $E=-n$, $F=\delta$ it takes the form:
\begin{multline}\label{eq:5F4-1bal1shift-8F7-VWP}
{}_{5}F_{4}\left.\!\left(\begin{matrix}
(A-1)/2-B,C,D,E,F+1
\\
(A+1)/2,A-B+1,C+D+E-A,F
\end{matrix}\right.\right)
\\
=\frac{\Gamma(1+A-C-E)\Gamma(1+A-D-E)\Gamma(1+A-C-D)\Gamma(1+A)}{\Gamma(1+A-C)\Gamma(1+A-D)\Gamma(1+A-E)\Gamma(1+A-C-D-E)}
\\
\times{}_{8}F_{7}\left.\!\left(\begin{matrix}A,1+A/2,B,C,D,E,A/2-\sigma+1,A/2+\sigma+1
\\
A/2,1+A-B,1+A-C,1+A-D,1+A-E,A/2-\sigma,A/2+\sigma
\end{matrix}\right.\right),
\end{multline}
where $-E\in\N$, and
$$
\sigma^2=\frac{A(A-2F)(A-2B-1)-2F(A-2B)}{4(A-2B-2F-1)}.
$$
The formula remains true for non-integer $E$ provided that both sides converge. Note also that we can regard $\sigma$ on the right hand side as an arbitrary number while $F$ on the left hand side is then easily expressed in terms of $\sigma^2$.

%10
\item Combination of Wang-Rathie transformation \eqref{eq:quadWangRathie3.1}  with Dougall's summation formula \eqref{eq:Dougallk} leads to a transformation of
general Saalsch\"{u}tzian ${}_5F_{4}$ with one unit shift to a particular very well-poised ${}_9F_{8}$ with two unit shifts. Renaming parameters according to $A=2\alpha-1$, $B=2\alpha-\beta-1$, $C=2\alpha-\gamma$, $D=2\alpha-b_1$,$E=2\alpha-b_2$, $F=-n$, $G=\delta$, it takes the form:
\begin{multline}\label{eq:5F4-1bal1shift-9F8-VWP}
{}_{5}F_{4}\left.\!\!\left(\!\begin{matrix}
A-B-C,D,E,F,G+1
\\
1+A-B,1+A-C,D+E+F-A,G
\end{matrix}\right.\right)
\\
=\frac{\Gamma(A+1)\Gamma(1+A-D-F)\Gamma(1+A-E-F)\Gamma(1+A-D-E)}{\Gamma(1+A-F)\Gamma(1+A-D)\Gamma(1+A-E)\Gamma(1+A-D-E-F)}
\\
\times{}_{9}F_{8}\left.\!\!\left(\!\begin{matrix}
A,1+A/2,B,C,D,E,F,A/2-\omega+1,A/2+\omega+1
\\
A/2,1+A-B,1+A-C,1+A-D,1+A-E,1+A-F,A/2-\omega,A/2+\omega\end{matrix}\right.\right),
\end{multline}
where $-F\in\N$ and
$$
\omega^2=\frac{GBC}{A-B-C-G}+\frac{A^2}{4}.
$$
If we let $-F\to\infty$ over integers in \eqref{eq:5F4-1bal1shift-9F8-VWP} we obtain a relation for general non-terminating ${}_4F_3$ with one unit shift in terms of a very well-poised ${}_8F_{7}(-1)$ with two unit shifts
\begin{multline}\label{eq:4F3unitshift-8F7}
{}_{4}F_{3}\left.\!\!\left(\!\begin{matrix}
A-B-C,D,E,G+1
\\
1+A-B,1+A-C,G
\end{matrix}\right.\right)
=\frac{\Gamma(A+1)\Gamma(1+A-D-E)}{\Gamma(1+A-D)\Gamma(1+A-E)}
\\
\times{}_{8}F_{7}\left.\!\!\left(\!\begin{matrix}
A,1+A/2,B,C,D,E,A/2-\omega+1,A/2+\omega+1
\\
A/2,1+A-B,1+A-C,1+A-D,1+A-E,A/2-\omega,A/2+\omega\end{matrix}\right|-1\right).
\end{multline}
Note also that we can regard $\omega$ on the right hand side as an arbitrary number while $G$ on the left hand side is easily expressed in terms of $\omega^2$.

%11
\item Combination of the Rakha-Rathie transformation \eqref{eq:quadRakhaRathie3.1}  with \eqref{eq:Bailey2} leads to a transformation of a special Saalsch\"{u}tzian ${}_5F_{4}$ with one unit shift to a particular second kind nearly-poised   ${}_6F_{5}$ with three unit shifts.  Renaming parameters according to $A=\alpha$, $B=\alpha-\beta-1/2$, $C=1+2\alpha-b$, $D=\delta$, we obtain
\begin{multline}\label{eq:5F4-1bal1shift-6F5-NP2-3shift}
{}_{5}F_{4}\left.\!\!\left(\!\begin{matrix}-n,A+1,A-B-1/2,C,D+1
\\
1+2A-B,\Delta(1+C-n,2),D\end{matrix}\right.\right)
\\
=-\frac{(1+2A-C)_n}{(C+n)(1-C)_{n-1}}{}_{6}F_{5}\left.\!\!\left(\!\begin{matrix}-n,2A,1+A,B,A-\sigma+1,A+\sigma+1
\\
A,1+2A-B,1+2A-C,A-\sigma,A+\sigma\end{matrix}\right.\right),
\end{multline}
where $n\in\N$, and
$$
\sigma^2=\frac{(A-B)(A(A-1)-D/2)+A(D-A)/2}{A-B-D-1/2}.
$$
Note also that we can regard $\sigma$ or $A-\sigma$ on the right hand side as an arbitrary number while $D$ on the left hand side is easy to express in terms of $\sigma^2$.

%12
\item Combination of Wang-Rathie transformation \eqref{eq:quadWangRathie3.1}  with \eqref{eq:Bailey2} leads to a transformation of a special
Saalsch\"{u}tzian ${}_6F_{5}$ with one unit shift to a particular second kind nearly-poised  ${}_7F_{6}$ with three unit shifts. Renaming parameters according to $A=\alpha$, $B=2\alpha-\beta-1$, $C=2\alpha-\gamma$, $D=2\alpha-b$, $E=\delta$, we get
\begin{multline}\label{eq:6F5-1bal1shift-7F6-NP2-3shift}
{}_{6}F_{5}\left.\!\!\left(\!\begin{matrix}-n,A,A+1/2,2A-B-C-1,D,E+1
\\
2A-B,2A-C,\Delta(1+D-n,2),E\end{matrix}\right.\right)
=-\frac{(2A-D)_n}{(D+n)(1-D)_{n-1}}\times
\\
{}_{7}F_{6}\left.\!\!\left(\!\begin{matrix}-n,2A-1,A+1/2,B,C,A-\omega+1/2,A+\omega+1/2
\\
A-1/2,2A-B,2A-C,2A-D,A-\omega-1/2,A+\omega-1/2\end{matrix}\right.\right),
\end{multline}
where $n\in\N$, and
$$
\omega^2=\frac{BCE}{2A-B-C-E-1}+\left(A-\frac{1}{2}\right)^2.
$$
Note also that we can regard $\omega$ on the right hand side as an arbitrary number while $E$ on the left hand side is easy to express in terms of $\omega^2$.

%13
\item  Combination of \eqref{eq:quad15.8.15}  with IPD summation \eqref{eq:IPD12} leads to the following identity:
\begin{multline}\label{eq:3plusIPD(-1)-4plusIPD}
F\left.\!\!\left(\!\begin{matrix}\alpha,\beta,d,\hh+\pp
\\
1-\beta+\alpha,e,\hh\end{matrix}\right|-1\right)
=\frac{\Gamma(e-d-\alpha)\Gamma(1+d+\alpha-e)\Gamma(e)}{\Gamma(e-\alpha)\Gamma(1+d+\alpha-e+p)(\hh)_{\pp}}
\\
\times F\left.\!\left(\!\begin{array}{l}\Delta(\alpha,2),d,1-e+\alpha
\\
1-\beta+\alpha,\Delta(1+d+\alpha-e+p,2)\end{array}\, \vdots\, Y_p(1,2)\, \right.\right),
\end{multline}
where $-d\in\N$ and $Y_p(1,2)$ is defined in \eqref{eq:Ypolynomial} with $\lambda=-\alpha$.

%14
\item  Combination of \eqref{eq:quad15.8.21}  with IPD summation \eqref{eq:IPD12} leads to the following transformation formula:
\begin{multline}\label{eq:3plusIPD-4plusIPD}
F\left.\!\left(\!\begin{matrix}\alpha,\beta,\Delta(d,2),\Delta(\hh+\pp,2)
\\
1-\beta+\alpha,\Delta(e,2),\Delta(\hh,2)\end{matrix}\right.\right)
=\frac{\Gamma(e-d-2\alpha)\Gamma(e)}{\Gamma(e-2\alpha)(1+d+2\alpha-e)_{p}(\hh)_{\pp}}
\\
\times F\left.\!\left(\!\begin{array}{l}\alpha,\alpha-\beta+1/2,d,1-e+2\alpha
\\
1-\beta+\alpha,\Delta(1+d+2\alpha-e+p,2)\end{array}\, \vdots \, Y_p(1,2)\, \right.\right),
\end{multline}
where $-d\in\N$ and $Y_p(1,2)$ is defined in \eqref{eq:Ypolynomial} with $\lambda=-2\alpha$.

%15
\item  Combination of Whipple's quadratic transformation \eqref{eq:quadWhipple}  with IPD summation  \eqref{eq:IPD12} leads to the following identity:
\begin{multline}\label{eq:4plusIPD-5plusIPD}
F\left.\!\left(\!\begin{matrix}\alpha,\beta,\delta,d,\hh+\pp
\\
1-\beta+\alpha,1-\delta+\alpha,e,\hh\end{matrix}\right.\right)
=\frac{\Gamma(e-d-\alpha)\Gamma(1+d+\alpha-e)\Gamma(e)}{\Gamma(e-\alpha)\Gamma(1+d+\alpha-e+p)(\hh)_{\pp}}
\\
\times F\left.\!\left(\!\!\begin{array}{l}\Delta(\alpha,2),1+\alpha-\beta-\delta,d,1-e+\alpha
\\
1-\beta+\alpha,1-\delta+\alpha,\Delta(1+d+\alpha-e+p,2)\end{array}\, \vdots\,  Y_p(1,2)\, \right.\right),
\end{multline}
where $-d\in\N$ and  $Y_p(1,2)$ is defined in \eqref{eq:Ypolynomial} with $\lambda=-\alpha$.

%16
\item Combination Miller-Paris transformation \eqref{eq:MP2013-6.3} with  Dougall's summation formula \eqref{eq:Dougallk} leads to a generalization of \eqref{eq:6F5(-1)-VWP-3F24-1Bal}.
Renaming parameters according to $A=\alpha$, $B=\beta$, $C=1+\alpha-b_1$, $D=1+\alpha-b_2$, $E=-n$, it takes the form
\begin{multline}\label{eq:doug-mill}
F\left.\!\left(\begin{matrix}A,1+A/2,B,C,D,E,\hat{\etta}+1
\\
A/2,1+A-B,1+A-C,1+A-D,1+A-E,\hat{\etta}\end{matrix}\:\right|-1\right)
=\frac{\Gamma(1+A-C-D-E)}{\Gamma(1+A)}
\\
\times\frac{\Gamma(1+A-C)\Gamma(1+A-D)\Gamma(1+A-E)}{\Gamma(1+A-C-D)\Gamma(1+A-D-E)\Gamma(1+A-C-E)}
F\left.\!\left(\begin{matrix}C,D,E,\f+\m
\\
1+A-B,C+D+E-A,\f\end{matrix}\right.\right),
\end{multline}
where $-E\in\N$ and $\hat{\etta}$ are the roots of $\hat{R}_{2m}(t;A,B,\f,\m)$ . This formula extends to general $E$ via Carlson's theorem. For $m=r=1$, we have $2\hat{\eta}_{1,2}=A\pm\sqrt{A^2-4fB}$.
In this case we get a connection between general ${}_4F_{3}$ with one unit shift and very well-poised ${}_8F_7(-1)$ with two unit shifts.
A similar connection is given by \eqref{eq:4F3unitshift-8F7}. These two transformations, however, are substantially different.

Taking $B=-1$, $-E\in\N$ in \eqref{eq:doug-mill}, we get a summation formula for Saalsch\"{u}tzian (or $1-$balanced) $_4F_3$ with one unit shift:
\begin{multline}\label{eq:4F3-1bal-1shift2}
_4F_3\left( \begin{matrix} C, D, E, f+1 \\ 2+A, C+D+E-A, f  \end{matrix}\right)
\\
=\frac{\Gamma(1+A)\Gamma(1+A-C-D)\Gamma(1+A-C-E)\Gamma(1+A-D-E)}{\Gamma(1+A-C)\Gamma(1+A-D)\Gamma(1+A-E)\Gamma(1+A-C-D-E)}
\\ \times\bigg(1-\frac{CDE(1+A-f)}{f(1+A-C)(1+A-D)(1+A-E)}\bigg).
\end{multline}

%17
\item Combination of  Maier's formula \eqref{eq:MaierTh3.1} with Dougall's summation formula \eqref{eq:Dougallk} leads to a generalization of Whipple's transformation \eqref{eq:7F6-VWP-4F3-1bal} with $k$-balanced ${}_4F_3$ on the right hand side. Renaming variables according to $A=\alpha$, $B=\beta$, $C=\delta$, $D=1+\alpha-b_1$, $E=1+\alpha-b_2$, $F=-n$, we can write this identity as follows:
\begin{multline}\label{q:VWPplusIPD-4F3-rbal}
F\left.\!\left(\begin{matrix}A,1+A/2,B,C,D,E,F,1-\rrho
\\
A/2,1+A-B,1+A-C,1+A-D,1+A-E,1+A-F,\rrho\end{matrix}\right.\right)
\\
=\frac{\Gamma(1+A-D)\Gamma(1+A-E)\Gamma(1+A-F)\Gamma(1+A-D-E-F)}{\Gamma(1+A)\Gamma(1+A-D-F)\Gamma(1+A-D-E)\Gamma(1+A-E-F)}
\\
\times{}_{4}F_{3}\left.\!\left(\!\begin{matrix}1+A-B-C-k,D,E,F
\\
1+A-B,1+A-C,D+E+F-A\end{matrix}\right.\right),
\end{multline}
where $-F\in\N$.  The formula is then extended to any values of parameters such that the left hand side converges while the right hand side terminates using Carlson's theorem.  The polynomial  $P_{2k}(t;A,B,C)$ is defined in  \eqref{eq:MaierQ2kpoly}.
For $k=1$, its roots are $2\rho_{1,2}=-A\pm\sqrt{A^2-4BC}$.
Using $B=-1$ and  $k=1$ in formula \eqref{q:VWPplusIPD-4F3-rbal} and assuming that  $-F\in\N$, we have a summation formula for $2-$balanced $_3F_2$:
	\begin{multline}\label{eq:4F3-2bal}
	_3F_2\left( \begin{matrix}  D,E,F \\ 2+A,  D+E+F-A \end{matrix}\right)
	\\
	=\frac{\Gamma(1+A)\Gamma(1+A-D-F)\Gamma(1+A-D-E)\Gamma(1+A-E-F)}{\Gamma(1+A-D)\Gamma(1+A-E)\Gamma(1+A-F)\Gamma(1+A-D-E-F)}
	\\ \times\bigg(1+\frac{DEF}{(1+A-D)(1+A-E)(1+A-F)}\bigg).
	\end{multline}	
This identity is equivalent to the formula  \cite[(3.1)]{KR2012} due to Kim and Rathie who extended  Saalsch\"{u}tzian summation formula for $_3F_2$ to $2-$balanced case.

%18
\item Combination of Choi and Rathie's quadratic transformation \eqref{eq:quadChRat2.1}  with IPD summation formula \eqref{eq:IPD12} leads to the following extension of the Karlsson-Minton summation theorem \eqref{eq:KarlssonMinton}:
\begin{multline}\label{eq:ChoiRathie2.1-IPD12}
F\left.\!\!\left(\!\begin{matrix}\alpha,\beta,d,\hh+\pp
\\
\beta+1,e,\hh\end{matrix}\right.\right)
=\frac{\Gamma(e-d-2\beta)\Gamma(e)}{\Gamma(e-2\beta)(1+d+2\beta-e)_{p}(\hh)_{\pp}}
\\
\times F\left.\!\left(\!\begin{array}{l}\beta,\beta-\alpha/2+1/2,\beta-\alpha/2+1,d,1-e+2\beta
\\
\beta+1,2\beta-\alpha+1,\Delta(1+d+2\beta-e+p,2)\end{array}\, \vdots\, Y_p(1,2)\, \right.\right),
\end{multline}
where $-d\in\N$, $p=p_1+\cdots+p_l$ and $Y_p(1,2)$ is defined in \eqref{eq:Ypolynomial} with $\lambda=-2\beta$.

%19
\item  Closely related to the previous transformation is the formula obtained by using \eqref{eq:quadChRat2.4} instead of \eqref{eq:quadChRat2.1} in combination with the IPD summation formula \eqref{eq:IPD12}:
\begin{multline}\label{eq:ChoiRathie2.4-IPD12}
F\left.\!\!\left(\!\begin{matrix}\alpha+1,2\alpha,\beta,d,\hh+\pp
\\
\alpha,\beta+1,e,\hh\end{matrix}\right.\right)
=\frac{\Gamma(e-d-2\beta)\Gamma(e)}{\Gamma(e-2\beta)(1+d+2\beta-e)_{p}(\hh)_{\pp}}
\\
\times F\left.\!\left(\!\begin{array}{l}\beta,\beta-\alpha,\beta-\alpha+1/2,d,1-e+2\beta
\\
\beta+1,2\beta-2\alpha+1,\Delta(1+d+2\beta-e+p,2)\end{array}\, \vdots\, Y_p(1,2)\, \right.\right),
\end{multline}
where $-d\in\N$, $p=p_1+\cdots+p_l$ and $Y_p(1,2)$ is defined in \eqref{eq:Ypolynomial} with $\lambda=-2\beta$.

%20
\item Combination of Rakha and Rathie's quadratic transformation \eqref{eq:quadRakhaRathie3.1}  with IPD summation formula \eqref{eq:IPD12} leads to the transformation:
\begin{multline}\label{eq:RakhaRathie3.1-IPD12}
F\left.\!\!\left(\begin{matrix}2\alpha,\alpha-\beta-1/2,1+\alpha-\sigma,1+\alpha+\sigma,d,\hh+\pp
\\
\alpha+\beta+3/2,\alpha-\sigma,\alpha+\sigma,e,\hh\end{matrix}\right.\right)
=\frac{\Gamma(e-d-2\alpha)\Gamma(e)}{\Gamma(e-2\alpha)(1+d+2\alpha-e)_{p}(\hh)_{\pp}}
\\
\times F\left.\!\left(\!\begin{array}{l}\alpha,\beta,\delta+1,d,1-e+2\alpha
\\
\alpha+\beta+3/2,\delta,\Delta(1+d+2\alpha-e+p,2)\end{array}\, \vdots\, Y_p(1,2)\, \right.\right),
\end{multline}
where $-\alpha\in\N$ and/or $-d\in\N$, $p=p_1+\cdots+p_l$,
$$
\sigma^2=\frac{\alpha^2\beta-\alpha\beta\delta-\beta\delta/2-\delta/4}{\beta-\delta}
$$
and $Y_p(1,2)$ is defined in \eqref{eq:Ypolynomial} with $\lambda=-2\alpha$.

%21
\item Combination of Miller and Paris' quadratic transformation \eqref{eq:MP2013-6.3} and Bailey's summation \eqref{eq:Bailey2} gives a generalization of \eqref{eq:4F3-NP2-4F3}:
\begin{multline}\label{eq:MP-Bailey}
F\left.\!\!\left(\!\begin{matrix}\alpha,1+\alpha/2,\beta,-n,\hat{\etta}+1
\\
\alpha/2,1-\beta+\alpha,\delta,\hat{\etta}\end{matrix}\right|-1\right)
=\frac{(\delta-\alpha-n-1)(\delta-\alpha)_{n-1}}{(\delta)_{n}}
\\
\times F\left.\!\!\left(\!\begin{matrix}1+\alpha/2,(1+\alpha)/2,1-\delta+\alpha,-n,\f+\m
\\
1-\beta+\alpha,(2-\delta+\alpha-n)/2,(3-\delta+\alpha-n)/2,\f\end{matrix}\right.\right),
\end{multline}
where $n\in\N$ and $\hat{\etta}$ are the roots of the polynomial $\hat{R}_{2m}(t;\alpha,\beta,\f,\m)$ defined in \eqref{eq:hatR2m}.
For $m=r=1$, these roots are $2\hat{\eta}_{1,2}=\alpha\pm\sqrt{\alpha^2-4f\beta}$.

%22
\item Combination of Maier's transformation \eqref{eq:MaierTh3.1} and Bailey's summation \eqref{eq:Bailey2} gives a generalization of Bailey's formula \eqref{eq:5F4-NP2-5F4-1bal}:
\begin{multline}\label{eq:MaierTh3.1-Bailey2}
F\left.\!\left(\!\begin{array}{l}\alpha,1+\alpha/2,\beta,\delta,-n
\\
\alpha/2,1+\alpha-\beta,1+\alpha-\delta,\gamma\end{array}\, \vdots\, P_{2k}\, \right.\right)
=\frac{(\gamma-\alpha-n-1)(\gamma-\alpha)_{n-1}}{(\gamma)_{n}}
\\
\times{}_{5}F_{4}\left.\!\!\left(\!\begin{matrix}1+\alpha/2,(1+\alpha)/2,\alpha-\beta-\delta-k+1,-n,1-\gamma+\alpha
\\
1+\alpha-\beta,1+\alpha-\delta,(2-\gamma+\alpha-n)/2,(3-\gamma+\alpha-n)/2\end{matrix}\right.\right),
\end{multline}
where  $n,k\in\N$ and the polynomial $P_{2k}(t;\alpha,\beta,\delta)$ is given in \eqref{eq:MaierQ2kpoly}.
For $k=1$ its roots are  $2\rho_{1,2}=-\alpha\pm\sqrt{\alpha^2-4\beta\delta}$.
The function ${}_{5}F_{4}$ on the RHS is $k+1$-balanced.

Setting $\beta=-1$ and $k=1$ in formula \eqref{eq:MaierTh3.1-Bailey2} we obtain  
$$
_4F_3\left( \begin{matrix} 1+\alpha/2, (1+\alpha)/2,  -n, 1-\gamma+\alpha\\ 2+\alpha,  \Delta(2-\gamma+\alpha-n,2) \end{matrix}\right)=\frac{(\gamma)_n}{(\gamma-\alpha-n-1)(\gamma-\alpha)_{n-1}}\times\bigg(1-\frac{n}{\gamma}\bigg).
$$	
The function on the left hand side is $2-$balanced $_4F_3$.

%23
\item Further generalization of the above transformation is obtained by using \eqref{eq:MaierTh3.4} instead of \eqref{eq:MaierTh3.1} and \eqref{eq:Bailey2} to sum the generalized hypergeometric function on the RHS of \eqref{eq:master1}:
\begin{multline}\label{eq:MaierTh3.4-Bailey2}
F\left.\!\left(\!\begin{array}{l}\alpha,1+\alpha/2,\beta,\delta,-n
\\
\alpha/2,1+\alpha-\beta,1+\alpha-\delta,\lambda\end{array}\, \vdots\, \hat{P}_{2k}\, \right.\right)
=\frac{(\lambda-\alpha-n-1)(\lambda-\alpha)_{n-1}}{(\lambda)_{n}}
\\
\times{}_{6}F_{5}\left.\!\!\left(\!\begin{matrix}1+\alpha/2,(1+\alpha)/2,\alpha-\beta-\delta-k+1,-n,1-\lambda+\alpha,\gamma+k
\\
1+\alpha-\beta,1+\alpha-\delta,(2-\lambda+\alpha-n)/2,(3-\lambda+\alpha-n)/2,\gamma\end{matrix}\right.\right),
\end{multline}
where $n,k\in\N$ and the polynomial $\hat{P}_{2k}(t;\alpha,\beta,\delta,\gamma)$ is defined in \eqref{eq:MaierhatQ2kpoly}.
The function ${}_{6}F_{5}$ on the RHS is Saalsch\"{u}tzian.

Setting $\beta=-1$ and $k=1$ in formula \eqref{eq:MaierTh3.4-Bailey2}, we get a summation formula for a particular Saalsch\"{u}tzian (or $1-$balanced) $_4F_3$ with one unit shift:
	\begin{multline}\label{eq:4F3-1bal-1shift3}
	_4F_3\left( \begin{matrix} 1+\alpha/2, (1+\alpha)/2, -n, 1-\lambda+\alpha, \gamma+1 \\ 2+\alpha,  \Delta(2-\lambda+\alpha-n,2), \gamma \end{matrix}\right)
	\\
	 =\frac{(\lambda)_n}{(\lambda-\alpha-n-1)(\lambda-\alpha)_{n-1}}\bigg(1+\frac{n(1+\alpha-\gamma)}{\lambda\gamma}\bigg),
	\end{multline}	
where $n\in\N$.
	
%24
\item Combination of  Maier's transformation \eqref{eq:MaierTh3.4} with \eqref{eq:Dougallk} leads to a generalization
of Whipple's transformation \eqref{eq:7F6-VWP-4F3-1bal}.  Renaming parameters according to $A=\alpha$, $B=\beta$, $C=\delta$, $D=1+\alpha-b_1$,  $E=1+\alpha-b_2$, $F=-n$, it can be written as
\begin{multline}\label{eq:MaierTh3.4-Dougallk}
F\left.\!\left(\!\begin{array}{l}A,1+A/2,B,C,D,E,F
\\
A/2,1+A-B,1+A-C,1+A-D,1+A-E,1+A-F\end{array}\, \vdots\, \hat{P}_{2k}\, \right.\right)
\\
=\frac{\Gamma(1+A-D)\Gamma(1+A-E)\Gamma(1+A-F)\Gamma(1+A-D-E-F)}{\Gamma(1+A)\Gamma(1+A-D-F)\Gamma(1+A-D-E)\Gamma(1+A-E-F)}
\\
\times{}_{5}F_{4}\left.\!\!\left(\!\begin{matrix}1+A-B-C,D,E,F,G+k
\\
1+A-B,1+A-C,D+E+F-A,G\end{matrix}\right.\right),
\end{multline}
where $-F,k\in\N$.  The formula is then extended to non-integer values of $F$  making both side convergent using Carlson's theorem.  The function on the right hand side is general Saalsch\"{u}tzian  ${}_{5}F_{4}$ with one integral shift.  The polynomial $\hat{P}_{2k}(t;A,B,C,G)$ is defined in \eqref{eq:MaierhatQ2kpoly}.
 For $k=1$ its roots are given by
$$
2\hat{\rho}_{1,2}=-A\pm\sqrt{A^2-4BCG/(B+C+G-A)}.
$$
 
%25
\item Combination of Wang and Rathie's formula \eqref{eq:quadWangRathie3.1} with IPD summation formula \eqref{eq:IPD12} yields:
\begin{multline}\label{eq:WangRathie3.1-IPD12}
F\left.\!\!\left(\begin{matrix}2\alpha-1,2\alpha-\beta-1,2\alpha-\gamma,1/2+\alpha-\omega,1/2+\alpha+\omega,d,\hh+\pp
\\
\beta+1,\gamma,\alpha-1/2-\omega,\alpha-1/2+\omega,e,\hh\end{matrix}\right.\right)
\\
=\frac{\Gamma(e-d+1-2\alpha)\Gamma(e)}{\Gamma(e-2\alpha+1)(d+2\alpha-e)_{p}(\hh)_{\pp}}
F\left.\!\left(\!\!\begin{array}{l}\alpha,\alpha-1/2,\beta+\gamma-2\alpha,\delta+1,d,2\alpha-e
\\
\beta+1,\gamma,\delta,\Delta(d+2\alpha-e+p,2)\end{array}\, \vdots\, Y_p(1,2)\, \right.\!\!\right),
\end{multline}
where $-d\in\N$, $p=p_1+\cdots+p_l$,
$$
\omega^2=\left(\alpha-1/2\right)^2-\frac{\delta(\gamma-2\alpha)(2\alpha-\beta-1)}{\beta+\gamma-2\alpha-\delta}.
$$
and $Y_p(1,2)$ is defined in \eqref{eq:Ypolynomial} with $\lambda=1-2\alpha$.

%26
\item Combination of Miller and Paris' quadratic transformation \eqref{eq:MP2013-6.3} and IPD summation formula \eqref{eq:IPD12} yields:
\begin{multline}\label{eq:MP-IPD12}
F\left.\!\!\left(\!\begin{matrix}\alpha,\beta,d,\hat{\etta}+1,\hh+\pp
\\
1-\beta+\alpha,e,\hat{\etta},\hh\end{matrix}\right|-1\right)
=\frac{\Gamma(e-\alpha-d)\Gamma(e)}{\Gamma(e-\alpha)(1+d-e+\alpha)_{p}(\hh)_{\pp}}
\\
\times F\left.\!\left(\!\!\begin{array}{l}\alpha/2,(\alpha+1)/2,d,1+\alpha-e,\f+\m
\\
1-\beta+\alpha,\Delta(1+d+\alpha-e+p,2),\f\end{array}\, \vdots\, Y_p(1,2)\, \right.\!\!\right),
\end{multline}
where $-d\in\N$ and $\hat{\etta}$ are the roots of the polynomial $\hat{R}_{2m}(t;\alpha,\beta,\f,\m)$ defined in \eqref{eq:hatR2m}, $p=p_1+\cdots+p_l$,
and $Y_p(1,2)$ is defined in \eqref{eq:Ypolynomial} with $\lambda=-\alpha$.
Formula \eqref{eq:MP-IPD12} extends to non-integer values of $d$.

%27
\item Combination of Maier's transformation \eqref{eq:MaierTh3.1} and and IPD summation formula \eqref{eq:IPD12} yields:
\begin{multline}\label{eq:MaierTh3.1-IPD12}
F\left.\!\left(\!\begin{array}{l}\alpha,\beta,\delta,d,\hh+\pp
\\
1+\alpha-\beta,1+\alpha-\delta,e,\hh\end{array}\, \vdots\, P_{2k}\, \right.\right)
=\frac{\Gamma(e-\alpha-d)\Gamma(e)}{\Gamma(e-\alpha)(1+d-e+\alpha)_{p}(\hh)_{\pp}}
\\
\times F\left.\!\left(\!\!\begin{array}{l}\alpha/2,(\alpha+1)/2,d,1+\alpha-\beta-\delta-k,1+\alpha-e
\\
1-\beta+\alpha,1-\delta+\alpha,\Delta(1+d+\alpha-e+p,2)\end{array}\, \vdots\, Y_p(1,2)\, \right.\!\!\right),
\end{multline}
where $-d,k\in\N$ and the polynomial $P_{2k}=P_{2k}(t;\alpha,\beta,\delta)$ is defined in \eqref{eq:MaierQ2kpoly}
and the polynomial $Y_p(1,2)$ is defined in \eqref{eq:Ypolynomial} with $\lambda=-\alpha$.  For $k=1$ the roots of $P_{2k}$ are
$$
2\rho_{1,2}=-\alpha\pm\sqrt{\alpha^2-4\beta\delta}.
$$
Formula \eqref{eq:MaierTh3.1-IPD12} extends to non-integer values of $d$.

%28
\item A generalization of the previous transformation is obtained by using \eqref{eq:MaierTh3.4} instead of
\eqref{eq:MaierTh3.1}. Combining  \eqref{eq:MaierTh3.4}  with  IPD summation formula \eqref{eq:IPD12} we get:
\begin{multline}\label{eq:MaierTh3.4-IPD12}
F\left.\!\left(\!\begin{array}{l}\alpha,\beta,\delta,d,\hh+\pp
\\
1+\alpha-\beta,1+\alpha-\delta,e,\hh\end{array}\, \vdots\, \hat{P}_{2k}\, \right.\right)
=\frac{\Gamma(e-\alpha-d)\Gamma(e)}{\Gamma(e-\alpha)(1+d-e+\alpha)_{p}(\hh)_{\pp}}
\\
\times F\left.\!\left(\!\!\begin{array}{l}\alpha/2,(\alpha+1)/2,d,1+\alpha-\beta-\delta-k,1+\alpha-e,\gamma+k
\\
1-\beta+\alpha,1-\delta+\alpha,\Delta(1+d+\alpha-e+p,2),\gamma\end{array}\, \vdots\, Y_p(1,2)\, \right.\!\!\right),
\end{multline}
where $-d\in\N$ and the polynomial $\hat{P}_{2k}(t;\alpha,\beta,\delta,\gamma)$ is defined in \eqref{eq:MaierhatQ2kpoly}, while
and $Y_p(1,2)$ is given by \eqref{eq:Ypolynomial} with $\lambda=-\alpha$.   For $k=1$ the roots of $\hat{P}_{2k}$ are
$$
2\hat{\rho}_{1,2}=-\alpha\pm\sqrt{\alpha^2-4\beta\delta\gamma/(\beta+\delta+\gamma-\alpha)},
$$
so that the left hand side can be written as a standard hypergeometric function using \eqref{eq:FP-roots}.
Formula \eqref{eq:MaierTh3.4-IPD12} extends to non-integer values of $d$.
Setting  $\pp=1$, $\hh=h$ to be scalars in \eqref{eq:MaierTh3.4-IPD12} and assuming that the parameters are subject to the relation $1+\alpha-\beta-\delta-k=-1$, we get the following exotic summation formula 
\begin{multline*}
F\left.\!\left(\!\begin{array}{l}\alpha,\beta,\delta,d,h+1
\\
1+\alpha-\beta,1+\alpha-\delta,e,h\end{array}\, \vdots\, \hat{P}_{2k}\, \right.\right)
=\frac{\Gamma(e-\alpha-d)\Gamma(e)}{\Gamma(e-d)\Gamma(e-\alpha)(1+d+\alpha-e)h}
\\
\times\left(\frac{d(\alpha)_2(e-\alpha-1)(\gamma+k)[(\alpha+2)(h-d)-(h+1)(e-d-1)]}{(1-\beta+\alpha)(1-\delta+\alpha)(2+d+\alpha-e)(3+d+\alpha-e)\gamma}+\alpha(h-d)-h(e-d-1)\right).
\end{multline*}	
For $k=1$ the function on the left hand side is ${}_7F_{6}$ with three unit shifts.

\end{enumerate}

\subsection{Case V: $(u,v)=(2,2)$}
Combination of Kummer's second quadratic transformation \eqref{eq:quad15.8.13} with IPD summation \eqref{eq:IPD22} leads to the transformation formula
\begin{multline}\label{eq:Kummer+IPD22}
F\!\left(\left.\!\begin{array}{l}\alpha,\beta,d,\hh+\pp\\2\beta,e,\hh\end{array}\right|2\right)
\\
=\frac{\Gamma(e-\alpha-d)\Gamma(e)}{\Gamma(e-\alpha)(1+d+\alpha-e)_{p}(\hh)_{\pp}}
F\!\left(\!\begin{array}{l}\Delta(\alpha,2),\Delta(d,2)\\\beta+1/2,\Delta(1+d+\alpha+p-e,2)\end{array}\, \vdots\, Y_p(2,2)\, \right),
\end{multline}
where $\alpha$ or $d$ is a negative integer, $Y_p(2,2)$ is defined in \eqref{eq:Ypolynomial} with $\lambda=-\alpha$.
In particular, for $p=1$, according to \eqref{eq:Y1} and \eqref{eq:Y1root}, we get
\begin{multline*}
{}_{4}F_{3}\!\left(\left.\!\!\begin{array}{l}\alpha,\beta,d,h+1\\2\beta,e,h\end{array}\right|2\!\right)
\\
=\!\frac{\Gamma(e-\alpha-d)\Gamma(e)(\alpha-e-d-1-d\alpha/h)}{\Gamma(e-\alpha)\Gamma(e-d)(1+d+\alpha-e)}
{}_{5}F_{4}\!\left(\!\begin{array}{l}\Delta(\alpha,2),\Delta(d,2),\xi_*+1\\\beta+1/2,\Delta(2+d+\alpha-e,2),\xi_*\end{array}\!\!\right),
\end{multline*}
where
$$
\xi_*=-\frac{\alpha(d-h)+h(e-d-1)}{2h-2(e-1)}
$$
is the negated root of  $Y_1(2,2)$.

Combination of the Miller-Paris quadratic transformation \eqref{eq:MP2013-6.1} with IPD summation \eqref{eq:IPD22} leads to
\begin{multline}\label{eq:MP2013-6.1+IPD22}
F\!\left(\left.\!\begin{array}{l}\alpha,\beta-m,\etta+1,d,\hh+\pp\\2\beta,\etta,e,\hh\end{array}\right|2\right)
\\
=\frac{\Gamma(e-\alpha-d)\Gamma(e)}{\Gamma(e-\alpha)(1+d+\alpha-e)_{p}(\hh)_{\pp}}
F\!\left(\!\begin{array}{l}\Delta(\alpha,2),\Delta(d,2),\f+\m\\\beta+1/2,\Delta(1+d+\alpha+p-e,2),\f\end{array}\, \vdots\, Y_p(2,2)\, \right),
\end{multline}
where $\alpha$ or $d$ is a negative integer, $\etta$ is the vector of roots of the polynomial $R_{2m}$ defined in \eqref{eq:R2m} and $Y_p(2,2)$ is defined in \eqref{eq:Ypolynomial} with $\lambda=-\alpha$.  For $m=1$ the roots of $R_2$ are given by:
$$
\eta_{1,2}=2f-1/2\pm2\sqrt{(f+1/4)^2-f\beta}.
$$

\subsection{Case VI: $(u,v)=(2,1)$}
Combination of \eqref{eq:quad15.8.14} with \eqref{eq:IPD21} leads to the transformation formula
\begin{multline}\label{eq:case6+IPD21}
F\!\left(\!\begin{array}{l}\alpha,\beta,d,\hh+\pp\\2\beta,e,\hh\end{array}\right)
\\
=\frac{\Gamma(e-\alpha/2-d)\Gamma(e)}{\Gamma(e-\alpha/2)(1+d+\alpha/2-e)_{p}(\hh)_{\pp}}
F\!\left(\!\begin{array}{l}\alpha/2,\beta-\alpha/2,\Delta(d,2)
\\
\beta+1/2,e-\alpha/2,1+d+\alpha/2+p-e\end{array}\, \vdots\, Y_p(2,1)\, \right),
\end{multline}
where $Y_p(2,1)$ is defined in \eqref{eq:Ypolynomial} with $\lambda=-\alpha/2$.
In particular, for $p=1$, according to \eqref{eq:Y1} and \eqref{eq:Y1root}, we get
\begin{multline*}
{}_{4}F_{3}\!\left(\!\!\begin{array}{l}\alpha,\beta,d,h+1\\2\beta,e,h\end{array}\!\!\right)
\\
=\!\frac{\Gamma(e-\alpha/2-d)\Gamma(e)((h-d)\alpha/2-(e-d-1)h)}{\Gamma(e-\alpha/2)\Gamma(e-d)(1+d+\alpha/2-e)h}
{}_{5}F_{4}\!\left(\!\begin{array}{l}\alpha/2,\beta-\alpha/2,\Delta(d,2),\xi_*+1\\
\beta+1/2,e-\alpha/2,2+d+\alpha/2-e,\xi_*\end{array}\!\!\right),
\end{multline*}
where the negated root $\xi_*=-\xi$ of $Y_1(2,1)$ is given by
$$
\xi_*=-\frac{(d-h)\alpha/2+(e-d-1)h}{h+d-2e+2}.
$$
Both  above  formulas remain valid for any parameters making both sides convergent, which can be justified using Carlson's theorem.

%%%%%%%%%%%%%%%%%%%%%%%%%%%%%%%%%%%%%%%%%%

\section{Concluding remarks}
In this paper we derived over forty transformation formulas for the generalized hypergeometric function evaluated at a fixed argument, all of them presented in Section~3.  Most of these identities are new. We also included several known ones to demonstrate the power of our approach.  For each transformation presented here we have conducted a thorough search of the literature to verify whether it is a guise of a known result. In a few cases when such connection was found we provided the corresponding reference and explanation.  We further presented several new summation formulas obtainable from our transformation identities. The idea behind our method is rather simple and generalizes naturally the beta integral method explained in detail in \cite{KrRao2003}: we integrate a known transformation formula against the density expressed in terms of Meijer's $G$ function and apply a known summation formula to the resulting series.   We think that some of the formulas presented in this work may serve as a generating relations for certain groups of hypergeometric transformations. In particular, formulas \eqref{eq:5F4-2unitshifts} and \eqref{eq:5F4Saal} may clearly serve as such generators. A group-theoretic study of the resulting family of transformations is one possible direction of further research in analogy with the study undertaken by us in \cite{KPMath2020} the staring point of which is formula \eqref{eq:4F3generate}.  The groups of this sort play an important role in mathematical physics. In particular, they constitute the key ingredient of a succinct description of the symmetries of Clebsh-Gordon's and Wigner's $3-j$, $6-j$ and $9-j$ coefficients from the angular momentum theory \cite{KrRao2004,RDN,Rao,RaoBook}. The summation formulas for the generalized hypergeometric function with integral parameter differences (IPD), on the other hand, appear in calculation of several integrals in high energy field theories and statistical physics \cite{ShpotSrivastava}. A further recent and  application of IPD-type summation formula is in the area of multiple orthogonal polynomials, see \cite{LL2020}. Another possible research direction motivated by the present investigation are summation formulas for the hypergeometric functions with non-linearly restricted parameters, a rather striking example of which is formula  \eqref{eq:4F3truntedsum}. Further specialization of parameters in some of the transformations presented here will probably lead to new summation formulas with non-linearly restricted parameters.   

We believe that techniques presented in this paper may have further applications and potential extensions. In particular, they may be applied to new transformations formulas for the generalized hypergeometric functions like \cite[formula (6)]{ChenChu}. Further, they can definitely be extended to $k$-hypergeometric functions \cite{DiazPariguan}.  Applications to cubic transformations, including certain new summation formulas can be found in our paper \cite{CKP-Lob}.


\begin{thebibliography}{99}
\bibitem{AAR} G.E.\:Andrews, R.\:Askey and R.\:Roy, Special
functions, Cambridge University Press, 1999.

\bibitem{AndrStan}G.E.\:Andrews and D.\:Stanton, Determinants in plane partition enumeration,
Europ. J. Combinatorics 19 (1998), 273--282.

\bibitem{Askey1994}R.\:Askey, A look at the Bateman project, in: W. Abikoff, J.S. Birman, K. Kuiken (Eds.),
The Mathematical Legacy of Wilhelm Magnus: Groups, Geometry, and Special Functions, in: Contemporary Mathematics,
vol.169, American Mathematical Society, Providence, RI, 1994, pp.29--43.

\bibitem{Bailey} W.N.\:Bailey, Generalized hypergeometric series,
Stecherthafner Service Agency, New York and London, 1964. Reprinted from:
Cambridge Tracts in Mathematics and Mathematical Physics, {\bf 32}, 1935.
 
\bibitem{Blaschke} P.\:Blaschke, Hypergeometric form of the fundamental theorem of calculus, 	preprint, 2018, arXiv:1808.04837.

\bibitem{CKP-Lob} M.A.C.\:Candezano, D.B.\:Karp and E.G.\:Prilepkina, Further applications of the $G$ function integral method, 
Lobachevskii Journal of Mathematics, Vol. 41 (2020), No. 5, pp. 747--762.


\bibitem{Chen}K.-W.\:Chen, Clausen's Series ${}_3F_2(1)$ with Integral Parameter Differences,
Symmetry 2021, 13(10), 1783.
 
\bibitem{ChenChu} X.\:Chen and W.\:Chu, Evaluation of nonterminating ${}_3F_2(1/4)$-series perturbed by
three integer parameters, Analysis and Mathematical Physics (2021), 11:67.

\bibitem{ChoChungYun}Y-K.\:Cho, S.-Y.\:Chung, H\:Yun, Rational Extension of the Newton Diagram for the Positivity of  ${}_1F_{2}$ Hypergeometric Functions and Askey–Szeg\"{o} Problem,  Constructive Approximation,  Volume 51 (2020), 49--72. 

\bibitem{Choi}J.\:Choi, Certain Applications of Generalized Kummer’s Summation Formulas for ${}_2F_1$, Symmetry 2021, 13(8), 1538.

\bibitem{ChoiRathie}J.\:Choi and A.K.\:Rathie, Quadratic transformations involving
hypergeometric functions of two and higher order, East Asian Math. J. 22 (2006), No. 1, pp. 71--77.

\bibitem{Chu2002} W.\:Chu, Inversion techniques and combinatorial identities: balanced hypergeometric series, Rocky Mountain Journal of Mathematics, Volume 32, Number 2 (2002), 561--587.

\bibitem{ChuWang2009}W.\:Chu,  X\:Wang, The modified Abel lemma on summation by parts and terminating hypergeometric series identities,
Integral Transforms and Special Functions, 20:2(2009), 93--118.


\bibitem{DiazPariguan} R.\:D\'{i}az and E.\:Pariguan, On hypergeometric functions and Pochhammer k-symbol,
 Divulgaciones Matem\'{a}ticas Volume 15 No. 2(2007),  179--192.


\bibitem{GesselStant}I.\:Gessel And D.\:Stanton, Strange Evaluations Of Hypergeometric Series, Siam J. Math. Anal.
Vol. 13, No. 2 (1982), 295--308.

\bibitem{GreeneKnuth}	D.H\:Greene, D.E.\:Knuth, Mathematics for the analysis of algorithms, 3rd edition, 	Progress in computer science and applied logic series,	Birkh\"{a}user,	1990

\bibitem{Haglund}J.\:Haglund, Rook theory and hypergeometric series, Adv Appl Math. 17(1996), 408--459.

\bibitem{KLJAT2017}D.\:Karp and J.L.\:L\'{o}pez, Representations of hypergeometric functions for arbitrary values of the parameters and their use, Journal of Approximation Theory, Volume 218(2017), 42--70.

\bibitem{KPSIGMA} D.\:Karp and E.\:Prilepkina, Hypergeometric differential equation and new identities for the coefficients of N{\o}rlund and B\"{u}hring, SIGMA 12 (2016), 052, 23 pages.

\bibitem{KPChapter2019} D.B.Karp and E.G.Prilepkina, Alternative approach to Miller-Paris transformations and their extensions,  pp.117-140 in Transmutation Operators and Applications (edited by V.V.Kravchenko and S.M.Sitnik), Springer Trends in Mathematics Series, Birkh\"{a}user, 2020.

\bibitem{KPResults2019}D.B.\:Karp and E.G.\:Prilepkina, Degenerate Miller-Paris transformations, Results in Mathematics, (2019), 74--94.

\bibitem{KPITSF2018} D.B.\:Karp and E.G.\:Prilepkina, Extensions of Karlsson--Minton summation theorem and some consequences of the first Miller--Paris transformation, Integral Transforms and Special Functions, Vol. 29, Issue 12, (2018), 955--970.


\bibitem{KPMath2020}D.B.\:Karp and E.G.\:Prilepkina, Transformations of the hypergeometric ${}_4F_3$ with one unit shift: a group theoretic study, Mathematics 2020, 8(11), 1966.

\bibitem{KR2012}Y.S.\:Kim, A.K.\:Rathie
A new proof of Saalsch\"{u}tz's theorem for the series ${}_3F_2(1)$ and its contiguous results with
applications, Commun. Korean Math. Soc. 27 (2012), No. 1, 129--135.

\bibitem{KRP2013}Y.S.\:Kim, A.K.\:Rathie and R.B.\:Paris,  An extension of Saalsch\"{u}tz's summation theorem for the series  ${}_{r+3}F_{r+2}$, Integral Transforms and Special Functions, Volume 24, Issue 11 (2013), 916--921.

\bibitem{KRP2014}Y.S.\:Kim, A.K.\:Rathie and R.B.\:Paris, On two Thomae-type transformations for hypergeometric series with integral parameter differences, Math. Commun. 19(2014), 111--118.

\bibitem{KRP2015}Y.S.\:Kim, A.K.\:Rathie and R.B.\:Paris, An alternative proof of the extended Saalsch\"{u}tz summation theorem for the ${}_{r+3}F_{r+2}(1)$ series with applications, Mathematical Methods in the Applied Sciences.
38:18(2015), 4891--4900.

\bibitem{Koepf}W.\:Koepf, Hypergeometric Summation. An Algorithmic Approach to Summation and Special Function Identities. Second Edition, Springer, 2014.


\bibitem{Koepf2007}W.\:Koepf, Bieberbach's conjecture, the de Branges and Weinstein functions and the Askey-Gasper inequality, Ramanujan Journal, Volume 13, Issue 1--3 (2007), 103--129.

\bibitem{KrRao2003}C.\:Krattenthaler, K.\:Srinivasa Rao, Automatic generation of hypergeometric identities by the beta integral method, Journal of Computational and Applied Mathematics 160 (2003), 159--173.

\bibitem{KrRao2004}C.\:Krattenthaler, K.\:Srinivasa Rao, On group theoretical aspects, hypergeometric transformations and symmetries of angular momentum coefficients, Symmetries in Science XI (2005), 355--375, Kluwer Acad. Publ., Dordrecht, 2004.

\bibitem{LL2020}H.\:Lima, A:.\:Loureiro, Multiple orthogonal polynomials associated with confluent hypergeometric functions
Journal of Approximation Theory 260 (2020) 105484.

\bibitem{LukeBook} Y.L.\:Luke, The special functions and their approximations. Volume 1. Academic Press, 1969.

\bibitem{Maier2019} R.S.\:Maier, Extensions of the classical transformations of the hypergeometric function ${}_3F_2$,
Advances in Applied Mathematics 105(2019), 25--47.

\bibitem{MP2013}A.R.\:Miller and R.B.\:Paris, Transformation Formulas For The Generalized Hypergeometric Function
with Integral Parameter Differences, Rocky Mountain Journal Of
Mathematics Volume 43, Number 1 (2013), 291--327.

\bibitem{Minton} B.M.\:Minton, Generalized hypergeometric functions at unit argument, J Math Phys, 12(1970), 1375--1376.

\bibitem{NIST}F.W.J.\:Olver, D.W.\:Lozier,  R.F.\:Boisvert and C.W.\:Clark (Eds.) NIST Handbook of Mathematical Functions, Cambridge
University Press, 2010.


\bibitem{PBM3} A.P.\:Prudnikov, Yu.A.\:Brychkov, O.I.\:Marichev, Integrals and Series: More Special Functions,  Volume~3, Gordon and Breach Science Publishers, New York, 1990.

\bibitem{RakhaRathie}M.A.\:Rakha and A.K.\:Rathie, Extensions of Euler type II transformation and Saalsch\"{u}tz's theorem,
Bull. Korean Math. Soc. 48 (2011), No. 1, 151--156. DOI:10.4134/BKMS.2011.48.1.151

\bibitem{Rao} K.S.\:Rao, Hypergeometric series and Quantum Theory of Angular Momentum, in {\emph Selected Topics in Special Functions}; Eds.: R.P.\:Agarwal, H.L.\:Manocha and K.\:Srinivasa Rao; Allied Publishers Ltd., 2001; pp. 93--134.

\bibitem{RaoBook} K.S.\:Rao and V.\: Lakshminarayanan, Generalized Hypergeometric Functions: transformations and group theoretical aspects, IOP Publishing, 2018.

\bibitem{RDN}K.S.\:Rao, H.D.\:Doebner, P.\:Nattermann, Generalized hypergeometric series and the symmetries of 3-$j$ and 6-$j$ coefficients, pp.381--403 in S.\:Kanemitsu and C.\:Jia(eds.), Number Theoretic Methods. Future Trends, Springer, 2002.

\bibitem{WangRathie}X.Wang, A.K.Rathie, Extension of quadratic transformation due to Whipple with an application,
Advances in difference equations , 2013:157.

\bibitem{Skwarczynski}M.\:Skwarczy\'{n}ski, De Branges theorem and generalized hypergeometric functions,  Bulletin de la Soci\'{e}t\'{e} des Sciences et des Lettres de \L\'{o}d\'{z}, Vol. LXI, no. 1(2011), 47--103.

\bibitem{Seaborn}J.B.\:Seaborn, Hypergeometric Functions and Their Applications, Springer, 1991.

\bibitem{ShpotSrivastava}  M.A.\:Shpot and H.M.\:Srivastava, The Clausenian hypergeometric function ${}_3F_2$ with unit argument and negative integral parameter differences,  Applied Mathematics and Computation, 2015, 259, 819--827.

\bibitem{Slater} L.J.\:Slater, Generalized Hypergeometric Functions, Cambridge University Press, 1966.

\end{thebibliography}
\end{document}